\documentclass[a4paper,12pt]{article}

\usepackage{amsmath}
\usepackage{amsfonts}
\usepackage{amssymb}
\usepackage{amsthm}
\usepackage{amsopn}
\usepackage{amscd}

\usepackage{pstricks}
\usepackage[utf8]{inputenc}
\usepackage{courier}
\usepackage{url}
\usepackage{color}
\usepackage{subfigure}
\usepackage{xypic}
\usepackage{enumitem}
\setenumerate{label=(\roman*)}

\usepackage[vcentering,dvips]{geometry}
\def\address#1#2{\begingroup
\noindent\parbox[t]{7.8cm}{%
\small{\scshape\ignorespaces#1}\par\vskip1ex
\noindent\small{\itshape E-mail address}%
\/: #2\par\vskip4ex}\hfill%
\endgroup}%

\author{Nathan Owen Ilten \& Robert Vollmert}
\title{Deformations of Rational T-Varieties}
\date{}
\newcommand{\floor}[1]{\lfloor#1\rfloor}
\newcommand{\pair}[2]{\langle#1,#2\rangle}

\newcommand{\emb}{\hookrightarrow}

\DeclareMathOperator{\tail}{tail}
\DeclareMathOperator{\spec}{Spec}
\DeclareMathOperator{\cone}{Cone}

\DeclareMathOperator{\tv}{TV}

\DeclareMathOperator{\conv}{Conv}

\DeclareMathOperator{\PP}{\mathbb{P}}
\DeclareMathOperator{\T}{\mathcal{T}}
\DeclareMathOperator{\aut}{Aut}

\DeclareMathOperator{\WDiv}{Div}
\DeclareMathOperator{\Div}{div}
\DeclareMathOperator{\id}{id}

\DeclareMathOperator{\loc}{Loc}
\DeclareMathOperator{\face}{face}
\DeclareMathOperator{\relint}{relint}
\DeclareMathOperator{\supp}{supp}

\newcommand{\CO}{\mathcal{O}}

\newcommand{\ZZ}{\mathbb{Z}}
\newcommand{\QQ}{\mathbb{Q}}
\newcommand{\NN}{\mathbb{N}}
\newcommand{\RR}{\mathbb{R}}
\newcommand{\Aff}{\mathbb{A}}
\newcommand{\CC}{\mathbb{C}}
\newcommand{\F}{\mathcal{F}}
\newcommand{\E}{\mathcal{E}}
\newcommand{\A}{\mathcal{A}}
\newcommand{\B}{\mathcal{B}}
\newcommand{\G}{\mathcal{G}}
\newcommand{\mcH}{\mathcal{H}}

\newcommand{\C}{\mathcal{C}}
\newcommand{\J}{\mathcal{J}}
\newcommand{\I}{\mathcal{I}}
\newcommand{\D}{\mathcal{D}}
\newcommand{\mcP}{\mathcal{P}}

\newcommand{\tot}{\mathrm{tot}}
\newcommand{\dfan}{\mathcal{S}}
\newcommand{\base}{B}

\newcommand{\DP}{\D_{\!P}}
\newcommand{\DQ}{\D_{\!Q}}
\newcommand{\aP}{a_{\!P}}
\newcommand{\Yt}{{Y^\tot}}
\newcommand{\Dt}{{\D^\tot}}

\newtheorem{lemma}{Lemma}[section]
\newtheorem{prop}[lemma]{Proposition}
\newtheorem{cor}[lemma]{Corollary}
\newtheorem{thm}[lemma]{Theorem}

\theoremstyle{definition}
\newtheorem{ex}[lemma]{Example}
\newtheorem{remark}[lemma]{Remark}
\newtheorem{defn}[lemma]{Definition}

\newcommand{\blowupfan}{%
 \psset{unit=0.5cm}
 \begin{pspicture}(-4,-3.2)(4,3.2)%
   \psgrid[gridwidth=0.3pt,griddots=5,subgriddiv=1,gridlabels=5pt](-3,-3)(3,3)

 \psset{linewidth=1pt}%

 \psline{-}(3,3)(-3,-3)%
 \psline{-}(-3,3)(3,-3)%
 \psline{-}(0,3)(0,-3)%
 \psline{-}(3,0)(-3,0)
 \psset{linecolor=gray,linestyle=dashed,linewidth=1.5pt}
 \psline{-}(3,1)(-3,1)
\psline{-}(3,-1)(-3,-1)
\end{pspicture}}

\newcommand{\blowupnull}{%
 \psset{unit=0.5cm}
 \begin{pspicture}(-4,-4)(4,2)%
\psgrid[gridwidth=0.3pt,griddots=5,subgriddiv=1,gridlabels=5pt](-3,-1)(3,1)
  \psset{linewidth=1pt}%
  \psline{<-]}(-3,0)(-1,0)
  \psline{[->}(1,0)(3,0)
  \psline{[-]}(1,0)(0,0)
  \psline{[-]}(-1,0)(0,0)
  \rput(2.2,1.6){$\mathcal{D}_0$}
  \rput(.6,1.6){$\mathcal{C}_0$}
  \rput(-.6,1.6){$\mathcal{B}_0$}
  \rput(-2.2,1.6){$\mathcal{A}_0$}

\end{pspicture}
}

\newcommand{\blowupinfty}{%
 \psset{unit=0.5cm}
 \begin{pspicture}(-4,-2)(4,4)%
\psgrid[gridwidth=0.3pt,griddots=5,subgriddiv=1,gridlabels=5pt](-3,-1)(3,1)
  \psset{linewidth=1pt}%
  \psline{<-]}(-3,0)(-1,0)
  \psline{[->}(1,0)(3,0)
  \psline{[-]}(1,0)(0,0)
  \psline{[-]}(-1,0)(0,0)
  \rput(2.2,1.6){$\mathcal{H}_\infty$}
  \rput(.6,1.6){$\mathcal{G}_\infty$}
  \rput(-.6,1.6){$\mathcal{F}_\infty$}
  \rput(-2.2,1.6){$\mathcal{E}_\infty$}

\end{pspicture}
}

\newcommand{\cotangfannull}{%
 \psset{unit=0.5cm}
 \begin{pspicture}(-4,-3.2)(4,3.2)%
   \psgrid[gridwidth=0.3pt,griddots=5,subgriddiv=1,gridlabels=5pt](-3,-3)(3,3)

 \psset{linewidth=1pt}%

 \psline{-}(0,3)(0,1)(2,3)%
 \psline{-}(2,3)(0,1)(0,0)(3,0)%
 \psline{-}(3,0)(0,0)(0,-1)(2,-3)%
 \psline{-}(0,-3)(0,-1)(2,-3)%
 \psline{-}(0,-3)(0,-1)(-3,-1)%
 \psline{-}(-3, -1)(0,-1)(0,0)(-3,3)%
 \psline{-}(-3, 3)(0,0)(0,1)(-1,3)%
\end{pspicture}}

\newcommand{\cotangfaninfty}{%
 \psset{unit=0.5cm}
 \begin{pspicture}(-4,-3.2)(4,3.2)%
   \psgrid[gridwidth=0.3pt,griddots=5,subgriddiv=1,gridlabels=5pt](-3,-3)(3,3)
  \psset{linewidth=1pt}%
  \psline{-}(-3,0)(3,0)
  \psline{-}(0,-3)(0,3)
  \psline{-}(-3,3)(0,0)(-1.5,3)
  \psline{-}(3,2)(1,0)(3,-2)
  \end{pspicture}}

\newcommand{\cotangfanone}{%
 \psset{unit=0.5cm}
 \begin{pspicture}(-4,-3.2)(4,3.2)%
\psgrid[gridwidth=0.3pt,griddots=5,subgriddiv=1,gridlabels=5pt](-3,-3)(3,3)
 \psset{linewidth=1pt}%
  \psline{-}(-3,3)(0,0)(3,0)
  \psline{-}(0,-3)(0,3)
  \psline{-}(3,3)(0,0)(3,-3)
  \psline{-}(-3,1)(-1,1)(-2,3)
\end{pspicture}}

\newcommand{\cotangfannulla}{%
 \psset{unit=0.5cm}
 \begin{pspicture}(-4,-3.2)(4,3.2)%
   \psgrid[gridwidth=0.3pt,griddots=5,subgriddiv=1,gridlabels=5pt](-3,-3)(3,3)

 \psset{linewidth=1pt}%

 \psline{-}(0,3)(0,1)(2,3)%
 \psline{-}(2,3)(0,1)(0,0)(3,0)%
 \psline{-}(3,0)(0,0)%
 \psline{-}(0,-3)(0,0)(3,-3)%
 \psline{-}(0,-3)(0,0)(-3,0)%
 \psline{-}(-3, 0)(0,0)(0,0)(0,3)%
 \psline{-}(-3, 3)(0,0)(0,1)(-1,3)%
\end{pspicture}}

\newcommand{\cotangfannullb}{%
 \psset{unit=0.5cm}
 \begin{pspicture}(-4,-3.2)(4,3.2)%
   \psgrid[gridwidth=0.3pt,griddots=5,subgriddiv=1,gridlabels=5pt](-3,-3)(3,3)

 \psset{linewidth=1pt}%

 \psline{-}(0,3)(0,0)(3,3)%
 \psline{-}(3,3)(0,0)(0,0)(3,0)%
 \psline{-}(3,0)(0,0)(0,-1)(2,-3)%
 \psline{-}(0,-3)(0,-1)(2,-3)%
 \psline{-}(0,-3)(0,-1)(-3,-1)%
 \psline{-}(-3, -1)(0,-1)(0,0)(-3,3)%
 \psline{-}(-3, 3)(0,0)(0,0)(-2,3)%
\end{pspicture}}

\newcommand{\cotangfaninftydegen}{%
 \psset{unit=0.5cm}
 \begin{pspicture}(-4,-3.2)(4,3.2)%
   \psgrid[gridwidth=0.3pt,griddots=5,subgriddiv=1,gridlabels=5pt](-3,-3)(3,3)
  \psset{linewidth=1pt}%
  \psline{-}(-3,3)(0,0)(3,0)
  \psline{-}(0,-3)(0,3)
  \psline{-}(-3,1)(-1,1)(-2,3)
  \psline{-}(3,2)(1,0)(3,-2)
  \end{pspicture}}

 \newcommand{\threefoldnull}{%
 \psset{unit=0.5cm}
 \begin{pspicture}(-4,-3.2)(4,3.2)%
   \psgrid[gridwidth=0.3pt,griddots=5,subgriddiv=1,gridlabels=5pt](-3,-3)(3,3)
  \psset{linewidth=1pt}%
  \psline{-}(-3,0)(3,0)
  \psline{-}(0,-3)(0,3)
  \psline{-}(-3,-3)(3,3)
  \psline{-}(-3,1)(0,1)(2,3)
  \psline{-}(-3,-1)(-1,-1)(-1,-3)
  \psline{-}(1,-3)(1,0)(3,2)
  \end{pspicture}}

  \newcommand{\threefoldinfty}{%
 \psset{unit=0.5cm}
 \begin{pspicture}(-4,-3.2)(4,3.2)%
   \psgrid[gridwidth=0.3pt,griddots=5,subgriddiv=1,gridlabels=5pt](-3,-3)(3,3)
  \psset{linewidth=1pt}%
  \psline{-}(-3,0)(3,0)
  \psline{-}(0,-3)(0,3)
  \psline{-}(-3,-3)(3,3)
  \end{pspicture}}

  \newcommand{\threefoldnulla}{%
 \psset{unit=0.5cm}
 \begin{pspicture}(-4,-3.2)(4,3.2)%
   \psgrid[gridwidth=0.3pt,griddots=5,subgriddiv=1,gridlabels=5pt](-3,-3)(3,3)
  \psset{linewidth=1pt}%
  \psline{-}(-3,0)(3,0)
  \psline{-}(0,-3)(0,3)
  \psline{-}(-3,-3)(3,3)
  \psline{-}(1,-3)(1,0)(3,2)
  \end{pspicture}}

\newcommand{\threefoldnullb}{%
 \psset{unit=0.5cm}
 \begin{pspicture}(-4,-3.2)(4,3.2)%
   \psgrid[gridwidth=0.3pt,griddots=5,subgriddiv=1,gridlabels=5pt](-3,-3)(3,3)
  \psset{linewidth=1pt}%
  \psline{-}(-3,0)(3,0)
  \psline{-}(0,-3)(0,3)
  \psline{-}(-3,-3)(3,3)
  \psline{-}(-3,1)(0,1)(2,3)
  \end{pspicture}}

  \newcommand{\threefoldnullc}{%
 \psset{unit=0.5cm}
 \begin{pspicture}(-4,-3.2)(4,3.2)%
   \psgrid[gridwidth=0.3pt,griddots=5,subgriddiv=1,gridlabels=5pt](-3,-3)(3,3)
  \psset{linewidth=1pt}%
  \psline{-}(-3,0)(3,0)
  \psline{-}(0,-3)(0,3)
  \psline{-}(-3,-3)(3,3)
  \psline{-}(-3,-1)(-1,-1)(-1,-3)
  \end{pspicture}}

\newcommand{\blowupnulla}{%
 \psset{unit=1cm}
 \begin{pspicture}(-3.1,-1.1)(3.1,2.1)%
\psgrid[gridwidth=0.3pt,griddots=5,subgriddiv=1,gridlabels=5pt](-3,-1)(3,1)
  \psset{linewidth=1pt}%
  \psline{<-]}(-3,0)(0,0)
  \psline{[->}(1,0)(3,0)
  \psline{[-]}(1,0)(0,0)
  \rput(2.2,1.6){$\D_0^0$}
  \rput(.5,1.6){$\C_0^0$}
  \rput(-.5,1.6){$\B_0^0$}
  \rput(-2.2,1.6){$\A_0^0$}
  \psset{linecolor=gray}\psline{->}(-.5,1.2)(-.05,.2)

\end{pspicture}
}
\newcommand{\blowupnullb}{%
 \psset{unit=1cm}
 \begin{pspicture}(-3.1,-1.1)(3.1,2.1)%
\psgrid[gridwidth=0.3pt,griddots=5,subgriddiv=1,gridlabels=5pt](-3,-1)(3,1)
  \psset{linewidth=1pt}%
  \psline{<-]}(-3,0)(-1,0)
  \psline{[->}(0,0)(3,0)
  \psline{[-]}(-1,0)(0,0)
  \rput(2.2,1.6){$\D_0^1$}
  \rput(.5,1.6){$\C_0^1$}
  \rput(-.5,1.6){$\B_0^1$}
  \rput(-2.2,1.6){$\A_0^1$}
\psset{linecolor=gray}\psline{->}(.5,1.2)(.05,.2)

\end{pspicture}
}
\newcommand{\blowupinftya}{%
 \psset{unit=1cm}
 \begin{pspicture}(-3.1,-1.1)(3.1,2.1)%
\psgrid[gridwidth=0.3pt,griddots=5,subgriddiv=1,gridlabels=5pt](-3,-1)(3,1)
  \psset{linewidth=1pt}%
  \psline{<-]}(-3,0)(0,0)
  \psline{[->}(1,0)(3,0)
  \psline{[-]}(1,0)(0,0)
  \rput(2.2,1.6){$\mcH_\infty^0$}
  \rput(.5,1.6){$\G_\infty^0$}
  \rput(-.5,1.6){$\F_\infty^0$}
  \rput(-2.2,1.6){$\E_\infty^0$}
\psset{linecolor=gray}\psline{->}(-.5,1.2)(-.05,.2)

\end{pspicture}
}
\newcommand{\blowupinftyb}{%
 \psset{unit=1cm}
 \begin{pspicture}(-3.1,-1.1)(3.1,2.1)%
\psgrid[gridwidth=0.3pt,griddots=5,subgriddiv=1,gridlabels=5pt](-3,-1)(3,1)
  \psset{linewidth=1pt}%
  \psline{<-]}(-3,0)(-1,0)
  \psline{[->}(0,0)(3,0)
  \psline{[-]}(-1,0)(0,0)
  \rput(2.2,1.6){$\mcH_\infty^1$}
  \rput(.5,1.6){$\G_\infty^1$}
  \rput(-.5,1.6){$\F_\infty^1$}
  \rput(-2.2,1.6){$\E_\infty^1$}
\psset{linecolor=gray}\psline{->}(.5,1.2)(.05,.2)

\end{pspicture}
}

\newcommand{\fanopoly}{%
 \psset{unit=0.5cm}
 \begin{pspicture}(-4,-0.2)(4,3.2)%
   \psgrid[gridwidth=0.3pt,griddots=5,subgriddiv=1,gridlabels=5pt](-3,0)(3,3)
 \psset{linewidth=1pt,hatchcolor=gray,fillstyle=vlines,hatchsep=.15,  linecolor=white}%
\pspolygon(-3,3)(-1,1)(1,1)(3,3)
 \psset{linecolor=black,linewidth=1pt}
   \psline{<->}(-3,3)(-1,1)(1,1)(3,3)
\end{pspicture}}

\newcommand{\fanopolya}{%
 \psset{unit=0.5cm}
 \begin{pspicture}(-4,-0.2)(4,3.2)%
   \psgrid[gridwidth=0.3pt,griddots=5,subgriddiv=1,gridlabels=5pt](-3,0)(3,3)
 \psset{linewidth=1pt,hatchcolor=gray,fillstyle=vlines,hatchsep=.15,  linecolor=white}%
\pspolygon(-3,3)(-1,1)(0,1)(2,3)
 \psset{linecolor=black,linewidth=1pt}
   \psline{<->}(-3,3)(-1,1)(0,1)(2,3)
\end{pspicture}}

\newcommand{\fanopolyb}{%
 \psset{unit=0.5cm}
 \begin{pspicture}(-4,0.2)(4,3.2)%
   \psgrid[gridwidth=0.3pt,griddots=5,subgriddiv=1,gridlabels=5pt](-3,0)(3,3)
 \psset{linewidth=1pt,hatchcolor=gray,fillstyle=vlines,hatchsep=.15,  linecolor=white}%
\pspolygon(-3,3)(0,0)(1,0)(3,2)(3,3)
 \psset{linecolor=black,linewidth=1pt}
   \psline{<->}(-3,3)(0,0)(1,0)(3,2)
\end{pspicture}}

\newcommand{\apictureforrobert}{%
 \psset{unit=0.5cm}
 \newgray{shadow}{.9}
 \begin{pspicture}(-2.2,-1.5)(5.2,6.7)%
	 \psset{linewidth=1pt}%
\pspolygon[linecolor=shadow](0,0)(5,0)(5,6.5)(0,6.5)
 \psset{linecolor=black,linewidth=1.2pt}
   \psline{<->}(0,-1)(5,-1)
   \psline{<->}(-1,6)(-1,0)
 \psset{linecolor=black,linewidth=1pt}
   \psline(0,4)(5,4)
   \psline(0,1)(5,1)
   \psline[linestyle=dashed](2.5,0)(2.5,6.5)
   \psline(0,1.5)(5,6.5)
	\psdots(-1,1)(-1,4)(2.5,-1)
	\rput(-2,1){$\infty$}
	\rput(-2,4){$0$}
	\rput(-2,6){$\mathbb{P}^1$}
	\rput(2.5,-1.6){$0$}
	\rput(6,-1){$\Aff^1$}
	\rput(7,4){$D^\tot(0,0)$}
	\rput(7,6.5){$D^\tot(0,1)$}
	\rput(7,1){$D^\tot(\infty,0)$}
\end{pspicture}}

\newcommand{\deltadecomp}{%
 \psset{unit=0.6cm}
	 \newgray{shadow}{.9}
 \begin{pspicture}(0,-.5)(2,1.5)%
 \psset{linewidth=1pt}%
\pspolygon[fillstyle=solid,fillcolor=shadow](0,0)(0,-.5)(1,-.5)(2,0)(2,1)(1.5,1.5)
\psdots(0,0)(1,0)(2,0)(1,1)(2,1)
\end{pspicture}}

\newcommand{\deltadecompa}{%
 \psset{unit=0.6cm}
	 \newgray{shadow}{.9}
 \begin{pspicture}(0,-.5)(1,1.5)%
 \psset{linewidth=1pt}%
\pspolygon[fillstyle=solid,fillcolor=shadow](0,0)(.5,.5)(1,0)
\psdots(0,0)(1,0)
\end{pspicture}}

\newcommand{\deltadecompb}{%
 \psset{unit=0.6cm}
	 \newgray{shadow}{.9}
 \begin{pspicture}(0,-.5)(1,1.5)%
 \psset{linewidth=1pt}%
\pspolygon[fillstyle=solid,fillcolor=shadow](0,0)(0,-.5)(1,0)(1,1)
\psdots(0,0)(1,0)(1,1)
\end{pspicture}}

\begin{document}
\maketitle
\begin{abstract}
We show how to construct certain homogeneous deformations for rational normal varieties with codimension one torus action. This can then be used to construct homogeneous deformations of any toric variety in arbitrary degree. For locally trivial deformations coming from this construction, we calculate the image of the Kodaira-Spencer map. We then show that for a smooth complete toric variety, our homogeneous deformations span the space of first-order deformations. 
\end{abstract}

\noindent Keywords: Toric varieties, deformation theory, $T$-varieties 
\\

\noindent MSC: Primary 14D15; Secondary 14M25.
\section*{Introduction}
There has been much progress made on understanding the deformation theory of toric varieties. The case of toric singularities has been studied extensively by K. Altmann, see for example \cite{MR1329519}, \cite{MR1452429}, and \cite{MR1798979}. There have been several results on deformations of non-affine toric varieties as well. In \cite{MR2092771} and \cite{MR2169828}, A. Mavlyutov constructed certain deformations of complete weak Fano toric varieties via, respectively, regluing an open cover with automorphisms, and representing one toric variety as a complete intersection inside of a larger toric variety. Furthermore, in~\cite{ilten08}, the first author constructed toric $\QQ$-Gorenstein deformations for partial resolutions of toric surface singularities.

More recently, the first author has provided a combinatorial description for the space of first-order deformations $T_X^1$ of a smooth complete toric variety $X$ in~\cite{ilten09a}. Additionally, in the case that $X$ is a surface, he constructed homogeneous deformations via Minkowksi decompositions of polyhedral subdivisions and showed that these deformations span $T_X^1$. In independent work, Mavlyutov presented a similar construction of certain homogeneous deformations in all dimensions~\cite{mavlyutov09a}.

The goal of this paper is to generalize the results of~\cite{ilten09a} in several directions. First of all, we construct multi-parameter deformations of  varieties of arbitrary dimension which are also not necessarily smooth. Secondly, we will look not only at deformations of toric varieties but also deformations of rational $T$-varieties of complexity one, that is, rational normal varieties admitting an effective codimension one torus action.    

Much as an $n$-dimensional toric variety can be described by an $n$-dimensional fan, an $n$-dimensional $T$-variety $X$ of complexity one can be described by a curve and some $n-1$-dimensional combinatorial data. We then construct a deformation of $X$ by somehow deforming the corresponding combinatorial data. In section~\ref{sec:tvar}, we give a short overview of the necessary theory of $T$-varieties. We then show how to construct homogeneous deformations of affine $T$-varieties in section~\ref{sec:affinecase}. Here we also describe the fibers of such deformations explicitly as $T$-varieties. Note that the deformation theory of affine $T$-varieties is being further developed by the second author in~\cite{vollmert09}.
 As a special case, we can of course consider toric varieties with an action by some subtorus. We describe this in detail in  section~\ref{sec:affinetoric} and show how to recover the deformations constructed by Altmann. In particular, we have a very natural description of toric deformations with non-negative degree, which are essential for constructing homogeneous deformations of complete toric varieties. 

In section~\ref{sec:fandecomp} we then show how to glue the deformations of affine $T$-varieties together to construct deformations of non-affine $T$-varieties. As in the affine case, we can also describe the fibers of such deformations explicitly as $T$-varieties. Restricting to the case of locally trivial deformations, we then calculate the Kodaira-Spencer map in section~\ref{sec:ks}. 

Of course, non-affine toric varieties provide again an example where our construction can be put to use. In section~\ref{sec:tc} we reformulate our Kodaira-Spencer calculation in nicer terms for this special case. For a smooth complete toric variety $X$, we then construct certain special deformations and show that they in fact span $T_X^1$. Thus, at least for smooth complete toric varieties, our  deformations provide a kind of skeleton of the versal deformation.

Our approach has some aspects in common with the independent work of Mavlyutov---both approaches construct deformations via Minkowski decomposition of some combinatorial data. However, an important difference can be found in the distinct ways in which we translate our combinatorial data into deformations. His construction relies on the homogeneous coordinate ring of a toric variety.  In contrast, our construction utilizes the language of $T$-varieties and polyhedral divisors.
\\
\\
\noindent\emph{Acknowledgements}: We would like to thank Hendrik S\"u\ss{} and Klaus Altmann for a number of helpful conversations. Thanks are also due to the anonymous referee for suggesting several improvements.


\section{T-Varieties}\label{sec:tvar}
 We recall several notions from~\cite{MR2426131}. As usual, let $N$ be a lattice with dual $M$ and let $N_\QQ$ and $M_\QQ$ be the associated $\QQ$ vector spaces.
For any polyhedron $\Delta\subset N_\QQ$, let $\tail(\Delta)$ denote its tailcone, that is, the cone of unbounded directions in $\Delta$. Thus, $\Delta$ can be written as the Minkowski sum of some bounded polyhedron and its tailcone. Now for $u \in \tail(\Delta)^\vee \cap M$, denote by $\face({\Delta},{u})$ the face of $\Delta$ upon which $u$ achieves its minimum.
We will also be considering the empty set as a polyhedron; any face of the empty set is itself the empty set.
A polyhedron $\Delta$ is \emph{nontrivial} if it is the empty set, or differs from its tailcone.
We will constantly assume that all polyhedra with which we deal contain no linear subspace of positive dimension.

Consider any smooth semiprojective variety $Y$ over $\CC$; recall that semiprojective means projective over some affine variety. By $\CC(Y)$ we denote its field of rational functions. For any $f\in\CC(Y)$, let $V(f)$ denote the divisor of zeros of $f$ in $Y$.

Let $\delta\subset N_\QQ$ be a pointed polyhedral cone.
 \begin{defn}A \emph{polyhedral divisor} on $Y$ with tail cone $\delta$ is a 
formal finite sum
$$\D = \sum_P \DP \otimes P,$$
where $P$ runs over all prime divisors on $Y$ and $\DP$ is a polyhedron with tailcone $\delta$.
Here, finite means that only finitely many coefficients are nontrivial. Note that the empty set is also allowed as a coefficient.
\end{defn}
We can evaluate a polyhedral divisor for every element $u \in \delta^\vee \cap M$ via
$$\D(u):=\sum_P \min_{v \in \DP} \langle v , u \rangle P$$
 in order to obtain a divisor on $Y$ with coefficients in $\QQ\cup\{\infty\}$, where the minimum over the empty set is defined to be $\infty$. We can also consider $\D(u)$ as a $\QQ$-divisor on $\loc \D:= Y \setminus \left( \bigcup_{\DP = \emptyset} P \right)$. By an abuse of notation, we shall write $H^0(Y,\D(u))$ for $H^0(\loc \D, \D(u)):=H^0(\loc \D, \CO(\D(u)))$.

\begin{defn}
A polyhedral divisor $\D$ is called {\em proper} if 
\begin{enumerate}
	\item $Y\setminus \loc \D$ is the support of a semiample divisor on $Y$;
	\item For all $u\in \delta^\vee\cap M$, $\D(u)$ is semiample on $\loc \D$;
	\item For all $u\in \relint \delta^\vee \cap M$,  $\D(u)$ is big on $\loc \D$.
\end{enumerate}
	\end{defn}

To a proper polyhedral divisor we associate an $M$-graded $k$-algebra and consequently
an affine scheme admitting a $T^N=N\otimes \CC^*$-action:
$$X(\D):= \spec \bigoplus_{u \in \delta^\vee \cap M} H^0(Y,\D(u))\cdot\chi^u.$$
This construction gives a normal variety of dimension $\dim Y+\dim N_\QQ$ together with a $T^N$-action.

We now wish to glue these affine schemes together; this requires some further definitions.
\begin{defn}\label{def:face}
Let $\D=\sum_P \DP \otimes P$, $\D'=\sum_P \DP' \otimes P$ be two proper polyhedral divisors on $Y$ with tail cones $\delta$ and $\delta'$.
\begin{enumerate}
\item We define their \emph{intersection} by $$\D \cap \D' := \sum_P (\DP \cap \DP') \otimes P.$$
\item We say  $\D' \subset \D$ if $\DP'\subset\DP$ for every prime divisor $P \in Y$.
\item For $y\in Y$ a not necessarily closed point, set $\D_y:=\sum _{y\in P} \DP$, where summation is via Minkowski addition.
\item \label{item:face} $\D'$ is a \emph{face} of $\D$ i.e. $\D' \prec \D$ if $\D'\subset \D$ and for each $y\in\loc(\D')$ there is a pair $(w_y,D_y)\in (\delta^\vee\cap M) \times  |\D(w_y)|$ such that $y\notin\supp(D_y)$, $\D_y'=\face(\D_y,w_y)$, and $\face(\D_v',w_y)=\face(\D_v,w_y)$ for all $v\in Y\setminus \supp(D_y)$.
\end{enumerate}
\end{defn}
If $\D' \subset \D$ then we have an inclusion
$$\bigoplus_{u \in (\delta')^\vee \cap M} H^0( Y,\D'(u))\cdot\chi^u
\supset \bigoplus_{u \in \delta^\vee \cap M} H^0( Y,\D(u))\cdot\chi^u$$ which corresponds to a dominant morphism $X(\D') \rightarrow X(\D)$. This is an open embedding exactly when $\D' \prec \D$.

\begin{defn}
A {\em divisorial fan} is a finite set $\dfan$ of proper polyhedral divisors such that for $\D,\D' \in \dfan$ we have $\D \succ \D' \cap \D \prec \D'$ with $\D' \cap \D$ also in $\dfan$.
For  a not necessarily closed point $y\in Y$, the \emph{slice} $\dfan_y$ is defined to be the set of all polyhedra $\D_y$ with $\D \in \dfan$.
\end{defn}

We may glue the affine varieties $X(\D)$ via
$$X(\D) \leftarrow X(\D \cap \D') \rightarrow X(\D').$$
This construction yields a normal scheme $X(\dfan)$ of dimension $\dim Y+\dim N_\QQ$ with a torus action by $T^N$. Note that all normal varieties with torus action can be constructed in this manner.

\begin{remark}
For any prime divisor $P$, the face condition ensures that $\dfan_P$ is in fact a polyhedral subdivision. If $Y$ is a curve, $X(\dfan)$ is complete if $Y$ is complete and $\dfan_P$ is a complete polyhedral subdivision for all points $P$.  In this case, we also say that $\dfan$ is complete.
\end{remark}

We will need the following lemma to construct deformation maps:
\begin{lemma}\label{lemma:regular}
Given a map $f\colon Y \to \base$ where $\base$ is affine, the composition of $f$ with the rational quotient map $X(\dfan) \dashrightarrow Y$ is regular.
\end{lemma}
\begin{proof}
The statement is local on $X(\dfan)$, hence we may assume $X = X(\D)$ is affine. Then $X$ is the affine contraction of $\spec_Y \bigoplus \CO(\D(u))\cdot \chi^u$ which maps to $\base$ regularly, inducing a regular map $X \to \base$.
\end{proof}

\begin{remark}
	If $Y$ is a smooth projective curve, some of the above definitions simplify. We first define the \emph{degree} of a polyhedral divisor by $$\deg \D := \sum_P \DP $$ where summation is via Minkowski addition.
	A polyhedral divisor $\D$ is then proper if and only if $\deg \D\subset \delta$, and for all $u\in\delta^\vee$ with $\min_{v\in \deg \D} \langle v,u\rangle= 0$ it follows that $u\notin \relint (\delta^\vee)$ and a multiple of $\D(u)$ is principal.
	Likewise, if $\D$ is proper, then $\D' \prec \D$ if and only if  $\DP'$ is a face of $\DP$ for every point $P \in Y$ and $\deg \D \cap \delta' = \deg \D'$, see \cite{ilten:09a} proposition 1.1.
\end{remark}

\begin{remark}\label{rem:downgrade}
	Let $N'$ be an $n$-dimensional lattice with dual $M'$, $\Sigma$ a fan in $N_\QQ'$, and $X=\tv(\Sigma)$ the associated toric variety. Choose some primitive $R\in M'$ and let $N=N'\cap R^\perp\subset N'$. Furthermore choose some cosection $s\colon N'\to N$. We can thus consider $X$ as a $T$-variety with codimension one torus action by the subtorus $T^{N}\subset T^{N'}$. In this case, $Y=\mathbb{P}^1$ and the divisorial fan $\dfan$ consists of polyhedral divisors $\D^\sigma$ for each $\sigma\in\Sigma$, where
	$$\D^\sigma=s(\sigma\cap [R=1])\otimes\{0\}+s(\sigma\cap [R=-1])\otimes\{\infty\}.$$
	Note that we define the set $[R=a]$ to be $\{v\in N'_\QQ| \langle v,R\rangle =a\}$, that is, the set of points in $N'_\QQ$ for which $R$ takes the value $a$.
	
	This \emph{downgrading} procedure can be generalized to consider the action of any subtorus $T$ of $T^{N'}$. If as above, $N$ is the cocharacter lattice of $T$ with cosection $s\colon N'\to N$, then $\tv(\Sigma)$ is described as a $T$-variety by a divisorial fan on $\tv(\Sigma_Y)$, where $\Sigma_Y$ is the coarsest common refinement of all cones $p(\sigma)$ for $\sigma\in\Sigma$. Here, $p:N'_\QQ\to (N'/N)_\QQ$ denotes the projection. The divisorial fan $\dfan$ once again consists of polyhedral divisors $\D^\sigma$ for each $\sigma\in\Sigma$, where now
	$$
	\D^\sigma=\sum_{\rho\in\Sigma_Y^{(1)}} s(\sigma\cap p^{-1}(\rho))\otimes D_\rho.
	$$
	$\Sigma^{(1)}$ denotes the set of rays of $\Sigma$, and by abuse of notation we denote a ray and its primitive lattice generator by the same symbol.
By $D_\rho$ we denote the torus invariant divisor corresponding to a ray $\rho$.
\end{remark}

\begin{figure}[htbp]
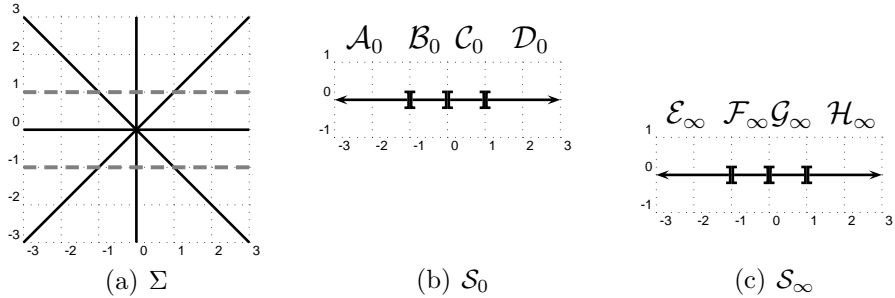

    \centering
    \subfigure[$\Sigma$]{\blowupfan}
    \subfigure[$\dfan_0$]{\blowupnull}
    \subfigure[$\dfan_\infty$]{\blowupinfty}
    \caption{Fan and divisorial fan for a toric surface}\label{fig:blowup}
   \end{figure}

\begin{ex}
	Consider the toric variety $X$ attained by blowing up $\mathbb{P}^1\times\mathbb{P}^1$ in all four fixpoints. The fan $\Sigma\subset N_\QQ'$ corresponding to this variety is pictured in figure~\ref{fig:blowup}(a). We choose $R=[0,1]$ and consider the sublattice $N=R^\perp$; the dashed gray lines in figure~\ref{fig:blowup}(a) mark  $[R=1]$ and $[R=-1]$. Choosing the cosection $s \colon N'\to N$ by $s(a,b)=(a)$ leads to the divisorial fan $\dfan$ pictured in figures~\ref{fig:blowup}(b) and~(c). The two-dimensional cones of $\Sigma$ correspond to the eight proper polyhedral divisors $\mathcal{A},\mathcal{B},\ldots,\mathcal{H}$. $X(\dfan)$ is thus our original variety $X$ considered with an action of the subtorus $T^N$.
\end{ex}

\begin{figure}[htbp]
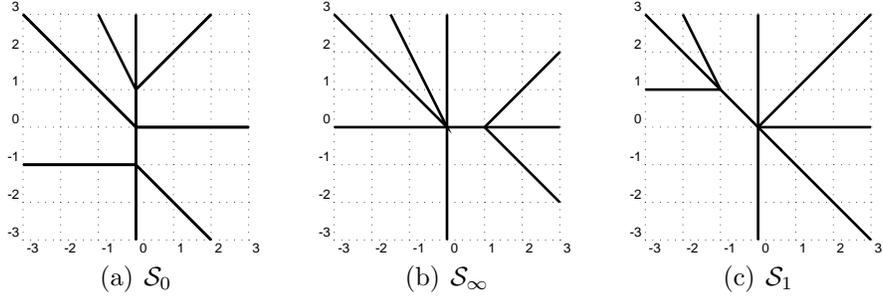

    \centering
    \subfigure[$\dfan_0$]{\cotangfannull}
    \subfigure[$\dfan_\infty$]{\cotangfaninfty}
    \subfigure[$\dfan_1$]{\cotangfanone}
    \caption{The divisorial fan for $\mathbb{P}(\Omega_{\mathcal{F}_1})$}\label{fig:cotanf1}
\end{figure}
\begin{ex}
	Let $\Omega_{\mathcal{F}_1}$ be the cotangent bundle of the first Hirzebruch surface. Then $X=\mathbb{P}(\Omega_{\mathcal{F}_1})$ is a $T$-variety over $\mathbb{P}^1$, see example 8.5 in~\cite{MR2426131}. The corresponding divisorial fan is pictured in figure~\ref{fig:cotanf1}, where two-dimensional polyhedra with common tailcone are coefficients for the same polyhedral divisor. 
\end{ex}

\section{Decompositions of Polyhedral Divisors}\label{sec:affinecase}
Let $Y$ be a smooth projective curve and let $\D$ be a proper polyhedral
divisor on $Y$ with $\delta = \tail(\D)$. We describe how to construct
deformations of $X=X(\D)$
with $T$ acting on the total space and preserving fibers. The construction is based on decomposing a coefficient
$\DP$ of $\D$ as a sum of polyhedra.

\begin{defn}
  An \emph{$r$-parameter Minkowski decomposition} of a 
  polyhedron $\Delta$ with tailcone $\delta$ is a decomposition
  $$
  \Delta = \Delta^0 + \dotsc + \Delta^r
  $$
  as Minkowski sum such that $\tail(\Delta^i)=\delta$ for $0\leq i \leq r$.
\end{defn}

\begin{defn}
  A Minkowski decomposition as above is said to be \emph{admissible} if
  it satisfies one of the following equivalent properties:
  \begin{enumerate}
  \item For each $u \in \delta^\vee \cap M$, at most one 
    $\face({\Delta^i},{u})$ has no lattice vertices.
  \item For each $u \in \delta^\vee \cap M$, at most one of the
    evaluations $\min\pair{\Delta^i}{u}$ is not an integer.
  \item For each vertex $v \in \Delta$, at most one of the corresponding
    vertices of the $\Delta^i$ is not a lattice point.
  \end{enumerate}
For example, 
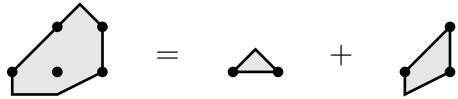
\begin{equation*}
	\begin{array}{c}\deltadecomp\end{array}\quad=\quad\begin{array}{c}\deltadecompa\end{array}\quad+\quad\begin{array}{c}\deltadecompb\end{array}
\end{equation*}
is an admissible one-parameter decomposition of a non-lattice polyhedron with tailcone $0$.
\end{defn}

\begin{remark}
We also consider Minkowski decompositions of $\Delta=\emptyset$, where we have fixed some cone $\delta$. In this case, an $r$-parameter Minkowski decomposition is defined identically to above; thus, at least one $\Delta^i$ must be the empty set.  Such a decomposition is admissible if for each $u\in\delta^\vee\cap M$, there is at most one $i$ with $\Delta^i\neq\emptyset$ and $\min\langle \Delta^i,u\rangle\notin\ZZ$.
\end{remark}

Let $\mathcal{P}\subset Y$ be a finite set of points in $Y$, including all those points $P$ with nontrivial coefficient $\DP$. Suppose now that for each $P\in \mathcal{P}$ we have Minkowski decompositions $\DP = \sum_{s=0}^{r_P} \DP^s$. We call such data a \emph{decomposition} of the polyhedral divisor $\D$; it is \emph{admissible} if each decomposition of the coefficients is admissible.  Let $r=\sum_{P\in\mcP} r_P$; this is finite since $\mcP$ is a finite set.

Consider some smooth affine variety $\base$ with special point $0$ cut out by a regular sequence $t_1,\ldots,t_{k}$. For $0\leq j \leq k$, let $\base_j$ be the subvariety cut out by $t_{j+1},\ldots,t_{k}$. Now consider some family $\gamma:Y^\tot \to \base$ with $Y^\tot$ smooth such that $Y_j^\tot:=V(t_{j+1},\ldots,t_{k})\subset Y^\tot$ is equal to $\gamma^{-1}(\base_j)$, and $\gamma^{-1}(0)=Y_0^\tot=Y$.
Furthermore, let $D^\tot(P,i)$ be a collection of pairwise different prime divisors on $Y^\tot$ intersecting the $Y_j^\tot$ properly, such
that $D^\tot(P,i)$ restricts to $P$ in $Y$. From this information, we then
define the polyhedral divisors
$$
\D^\tot = \sum_{P,i} \DP^i \otimes D^\tot(P,i).
$$
Note that since we required the $D^\tot(P,i)$ to restrict to $P$ in $Y$, we have
$\D^\tot_{|Y}=\D$, since  for each $P$, the
coefficients of the $D^\tot(P,i)$ sum up to $\DP$. In particularly,
$\D^\tot(u)_{|Y} = \D(u)$ for all $u$.

We assume for the moment that all $\D^\tot_{|Y_i^
\tot}$ are proper polyhedral divisors. Let $X^
\tot=X(\D^\tot)$.
By lemma~\ref{lemma:regular}, we get a map $\pi\colon X^\tot \to \base$.
We want the special fiber of $
\pi$ to be $X$, i.e. $\pi^{-1}(0)=X$.

\begin{prop}\label{prop:specialfiber}
The map of $T$-varieties $X \to X^\tot$ induced by 
$Y \emb Y^\tot$ embeds $X$ as the special fiber $\pi^{-1}(0)$ if, for each
$u \in \delta^\vee \cap M$, the following two conditions hold:
\begin{enumerate}
\item \label{cond:admissible}
$\floor{\D^\tot(u)}_{|Y} = \floor{\D^\tot(u)_{|Y}}$
\item \label{cond:surj}
With $D = \floor{\D^\tot(u)}$, the natural morphisms
$H^0(Y_i^\tot, D_{|Y_{i}^\tot}) \to H^0(Y_{i-1}^\tot, D_{|Y_{i-1}^\tot})$ are surjective for $1\leq i \leq k$.
\end{enumerate}
\end{prop}

\begin{proof}
The claim is equivalent to the exactness of 
\begin{equation*}	\begin{CD}	
0 @>>> I\cdot H^0(\Yt,\Dt(u))@>>>  H^0(\Yt,\Dt(u)) @>\nu>> H^0(Y,\D(u)) @>>> 0
\end{CD}
\end{equation*}
for each $u \in \delta^\vee \cap M$, where $I=\langle t_1,\ldots,t_k\rangle$. The  map $\nu$ arises as
follows (compare section~8 of~\cite{MR2426131}):
$$
\xymatrix{
H^0(Y^\tot,\D^\tot(u)) \ar@{=}[d] \ar[rr] & & H^0(Y,\D(u)) \ar@{=}[d] \\
H^0(Y^\tot,\floor{\D^\tot(u)}) \ar[r]^\varphi &
H^0(Y,\floor{\D^\tot(u)}_{|Y}) \ar@^{(->}[r]^\psi &
H^0(Y,\floor{\D(u)})
}
$$
Since $\Dt(u)_{|Y} = \D(u)$, surjectivity of $\psi$ follows from
condition~\ref{cond:admissible}.
Surjectivity of $\varphi$ follows from condition~\ref{cond:surj}. Thus, $\nu$ is surjective (and $X\to X^\tot$ is a closed embedding).

We now must check that the kernel of $\nu$ is correct; an easy calculation shows that it contains  $I\cdot H^0(\Yt,\Dt(u))$. Choose some open affine $U\subset Y^\tot$ such that $U\cap Y\neq \emptyset$ and $U$ is disjoint from the support of $D=\floor{\D^\tot(u)}$. Then we can expand the above sequence to 
\begin{equation*}	\begin{CD}	
	0 @>>> I\cdot H^0(\Yt,D)@>>>  H^0(\Yt,D) @>>> H^0(Y,D_{|Y}) @>>> 0\\
@. @VVV @VVV @VVV\\
0 @>>> I\cdot H^0(U,\CO_U)@>>>  H^0(U,\CO_U) @>>> H^0(Y\cap U,\CO_{Y\cap U}) @>>> 0
\end{CD}
\end{equation*}
where the vertical arrows are inclusions. If we can show that 
$ I\cdot H^0(U,\CO_U) \cap H^0(\Yt,D) = I\cdot H^0(\Yt,D)$, we are done by the exactness of the second row.

Assume that $k=1$ and take $s\in  I\cdot H^0(U,\CO_U) \cap H^0(\Yt,D) $; we can thus write $s=t_1 g $ for $g\in H^0(U,\CO_U) $. Furthermore,
$$
\Div (t_1 g)+D\geq 0.
$$
But $\Div (t_1 g)=\Div(t_1)+\Div(g)$ and the order of the components of $D$ along $\Div(t_1)=Y$ are zero, so $g\in H^0(\Yt,D)$. Thus $ I\cdot H^0(U,\CO_U) \cap H^0(\Yt,D) = I\cdot H^0(\Yt,D)$.

Assume that $k>1$. After slight adjustment, the above arguments show that 
\begin{equation}\label{eqn:indstep}	
	0 \to t_i\cdot H^0(Y_i^\tot,D_{|Y_i^\tot})\to  H^0(Y_i^\tot,D_{|Y_i^\tot}) \to H^0({Y_{i-1}^\tot},D_{|{Y_{i-1}^\tot}}) \to 0
\end{equation}
is exact for $1\leq i \leq k$. Now, consider also the sequence 
\begin{equation}\label{eqn:indass}	
	0 \to \langle t_1,\ldots,t_j\rangle \cdot H^0(Y_j^\tot,D_{|Y_j^\tot})\to  H^0(Y_j^\tot,D_{|Y_j^\tot}) \to H^0({Y,D_{|Y}}) \to 0
\end{equation}
and assume that this is exact for some $j=l$, $1\leq l < k$.
A straightforward diagram chase shows that the exactness of \eqref{eqn:indstep} for $i=l+1$  and exactness of \eqref{eqn:indass} for $j=l$ gives the exactness of \eqref{eqn:indass} for $j=l+1$.
Induction on $l$ completes the proof.
\end{proof}

Condition~\ref{cond:admissible} is where admissibility comes in to play:

\begin{lemma}
Suppose $D = \sum \aP^i D^\tot(P,i)$ is a $\QQ$-divisor on $Y^\tot$. Then
$\floor{nD}_{|Y} = \floor{(nD)_{|Y}}$ for all integers $n \ge 0$
if and only if, for each $P \in Y$, at most one of the coefficients
$\aP^i$ is not an integer.
\end{lemma}
\begin{proof}
Due to our choice of divisors $D^\tot(P,s)$, this follows from the
following fact: Let $p, q \in \QQ \setminus \ZZ$, $p, q \ge 0$. Then there
exists an integer $n \ge 0$ such that $\floor{np+nq} > \floor{np} + \floor{nq}$.
\end{proof}

\begin{cor}\label{cor:admissible}
Condition~\ref{cond:admissible} of proposition~\ref{prop:specialfiber}
holds for each $u \in \delta^\vee \cap M$ if and only if the Minkowski
decompositions underlying $\D^\tot$ are admissible.
\end{cor}

From now on, we assume that our base curve $Y$ is $\PP^1$. For each $P\in \mathcal{P}$, let $y_P\in \CC(Y)$ be a rational function with its sole zero at $P$. Let $t_{P,1},\ldots,t_{P,r_P}$ be coordinates on $\Aff^{r_P}$ for $P\in\mcP$, and set $t_{P,0}=0$.  Let $\base$ be any open affine neighborhood of the origin in $\prod_{P\in\mcP} \mathbb{A}^{r_P}$ such that a divisor on $\PP^1\times \base$ of the form $V(y_P-t_{P,i})$ doesn't intersect any divisor of the form  $V(y_Q-t_{Q,j})$ for $P\neq Q$ and $P,Q\in \mcP$.

\begin{figure}[htbp]
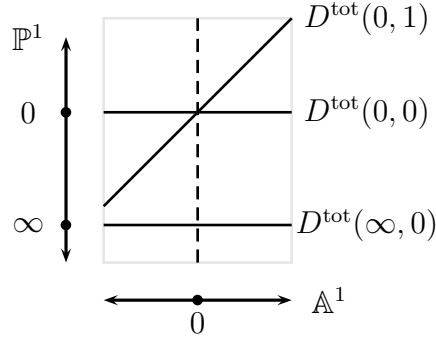

    \centering
    \apictureforrobert
    \caption{A family of prime divisors $D^\tot(P,i)$ on $\mathbb{P}^1\times\Aff^1$}\label{fig:p1a1}
\end{figure}

We now consider the trivial family $\Yt=\PP^1\times\base$.
As prime divisors, we then take $D^\tot(P,i) = V(y_P - t_{P,i})$; these clearly restrict as desired to $P$. As an example, such a family is pictured in figure~\ref{fig:p1a1} for $Y^\tot=\mathbb{P}^1\times \Aff^1$, with $r_0=1$ and $r_\infty=0$.

For a point $\lambda\in \base$, let $\D^{(\lambda)}$ be the restriction of $\D^\tot$ to $Y_\lambda^\tot$, the fiber over $\lambda$. Since $Y_\lambda^\tot\cong Y$, we can view $\D^{(\lambda)}$ as a polyhedral divisor on $Y$. In fact, we can describe $\D^{(\lambda)}$ explicitly. Say $\lambda$ is given by the equations $t_{P,i}=\lambda_{P,i}$, and set $\lambda_{P,0}=0$ for each $P\in\mathcal{P}$. For $0\leq i \leq r_P$, let $D^{(\lambda)}(P,i)$ be the divisor on $Y$ given by the vanishing of $y_P-\lambda_{P,i}$. Then the polyhedral divisor $\D^{(\lambda)}$ is given by
$$\D^{(\lambda)}=\sum_{\substack{P\in\mathcal{P}\\0\leq i\leq r_P}} \DP^{i}\otimes D^{(\lambda)}(P,i),$$
where the coefficients in front of prime divisors appearing multiple times are added via Minkowski sums. Note that prime divisors appear multiple times whenever $\lambda_{P,i}=\lambda_{P,j}$ for $i\neq j$.

\begin{lemma}\label{lemma:properness}
	$\D^\tot$ is a proper polyhedral divisor on $Y^\tot = \PP^1 \times \base$. Likewise, $\D^{(\lambda)}$ is a proper polyhedral divisor on $\mathbb{P}^1$. 
	\begin{proof}
		For any $u\in \delta^\vee \cap M$, consider the $\QQ$-divisor $\D^\tot(u)$. One
		easily checks that $l\cdot\D^\tot(u)\sim l\cdot\D(u)\times \base$ for some large $l\in\NN$, since for each
		$P\in Y$ and $i\leq r_P$ we have $V(y_P-t_{P,i})\sim V(y_P)$. Thus,
		$\D^\tot(u)$ is semiample or big exactly when $\D(u)$ is semiample
		or big, so the properness of $\D^\tot$ follows from the properness
		of $\D$. A similar argument holds for $\D^{(\lambda)}$.
	\end{proof}
\end{lemma}

Since $\D^\tot$ is proper it defines an affine $T$-variety $X^\tot=X(\D^\tot)$ and we have a natural map $\pi\colon X^\tot\to \base$. If the decomposition of $\D$ is admissible, this  is in fact a deformation:

\begin{thm}\label{thm:affinemainthm}
	If the decompositions of $\DP$ for each $P\in\mathcal{P}$ are admissible, the map $\pi\colon X(\D^\tot)\to \base$ gives a flat family with $\pi^{-1}(\lambda)\cong X(\D^{(\lambda)})$ for $\lambda\in \base$. In particular, $\pi^{-1}(0)=X(\D)=X$.
\end{thm}

\begin{remark}We call deformations of the above sort \emph{$T$-deformations}. We also say that they are \emph{homogeneous} since the torus $T$ acts on the total space, preserving fibers and extending the action of $T$ on the special fiber. 
\end{remark}
	The proof of the theorem relies on the following proposition:

\begin{prop}\label{prop:surjectivity}
	If the decomposition of $\D$ is admissible, the map $X \to X^\tot$ is a closed embedding given by the ideal generated by all $t_{P,i}$ for $P\in\mathcal{P}$, $1\leq i \leq r_P$.
\end{prop}

\begin{proof}
Let $D = \floor{\D^\tot(u)}$ for some $u \in \delta^\vee \cap M$.
We twist the short exact sequence for the embedding $Y_{i-1}^\tot \emb Y_i^\tot$
$$
0 \to \I_{Y_{i-1}^\tot} \to \CO_{Y_i^\tot} \to \CO_{Y_{i-1}^\tot} \to 0
$$
by the locally free sheaf $\CO_{Y^\tot}(D)_{|Y_i}$. Consider the associated
long exact sequence in cohomology
$$
H^0(Y_{i}^\tot,D_{|Y_{i}^\tot}) \to H^0(Y_{i-1}^\tot,D_{|Y_{i-1}^\tot}) \to H^1(Y_i^\tot,\I_{Y_{i-1}^\tot}(D)).
$$
Assume that $H^0(Y_{i-1}^\tot,D_{|Y_{i-1}^\tot})$ is not zero. 
We claim that $H^1(Y_i^\tot,\I_{Y_{i-1}^\tot}(D))$ vanishes, which proves the statement
by proposition~\ref{prop:specialfiber}.

Indeed, since $\I_{Y_{i-1}^\tot} =t_i\cdot \CO_{Y_i^\tot} $, we then have $H^1(Y_i^\tot,\I_{Y_{i-1}^\tot}(D))=H^1(Y_i^\tot,D_{|Y_{i}^\tot})$. But this disappears, since $D_{|Y_{i}^\tot}$ is an effective semiample divisor on the product of $\mathbb{P}^1$ with some subset of affine space.
\end{proof}

\begin{proof}[Proof of theorem \ref{thm:affinemainthm}]
	By choice of $\base \subset \Aff^r$, the admissible decomposition of $\D$ also induces admissible decompositions of $\D^{(\lambda)}$ which result in the same divisor $\D^\tot$ on $Y^\tot$. Thus, after coordinate change in $\base$ we can apply proposition~\ref{prop:surjectivity} and get that for any $\lambda\in \base$, $X(\D^\lambda)\cong \pi^{-1}(\lambda)$. Furthermore,  $X(\D^{(\lambda)}) \emb X^\tot$ is cut out by a regular sequence so $\pi$ is flat.
\end{proof}

\begin{remark}\label{rem:nonprim}
	Suppose a Minkowski summand $\DP^i$ is a multiple $k\Delta$
	of a lattice polyhedron $\Delta$.
	Then replacing $D^\tot(P,i) = V(y_P-t_{P,i})$
	by $D^\tot(kP,i) := V(y_P^k-t_{P,i})$ and $\DP^i$ by $\Delta$
	in $\D^\tot$ also gives a deformation of $X$, after changing $\base$ accordingly.
	Indeed,
	since $\Delta \otimes D^\tot(kP,i)$ restricts to
	$\Delta \otimes kP$, $\D^\tot$ restricts to $\D$ as before.
	The change doesn't affect the integrality considerations
	since $\Delta$ is a lattice polyhedron.
	The rest of the arguments carry through unchanged.
\end{remark}

We end this section with a corollary of theorem~\ref{thm:affinemainthm}:
\begin{cor}
	Let $\D$ be a proper polyhedral divisor on $\mathbb{P}^1$ with affine locus. Consider some admissible decomposition of $\D$. Then the general fiber of the corresponding deformation $\pi$ has exactly the analytic singularities $\tv(\cone(\DP^i\times\{1\}))$ for $P\in\mcP$ and $0\leq i \leq r_P$, where $\tv(\cone(\DP^i\times\{1\}))$ is the toric singularity corresponding to the cone over the polyhedron $\DP^i$.
	\begin{proof}
		This follows from the description of the general fiber from theorem~\ref{thm:affinemainthm} coupled with~\cite{liendo:10a}, theorem~5.3.
	\end{proof}
\end{cor}

\begin{remark}
	In section 4 of~\cite{ilten08}, explicit equations were used to calculate the singularities in the general fiber for toric deformations of cyclic quotient singularities. Combining this with the description of affine toric deformations in the following section, the above corollary provides a way of doing this without using the equations. Furthermore, the above corollary can be applied to see whether a toric deformation, or more generally, a $T$-deformation, is a smoothing. Note that if $\D$ has complete locus and $X(\D)$ is singular, no $T$-deformation can be a smoothing (see~\cite{liendo:10a}, proposition~5.1). 
\end{remark}

\section{Deformations of Affine Toric Varieties}\label{sec:affinetoric}

We turn now to the case that $X = \tv(\sigma)$ is a toric variety
with embedded torus $T'=T^{N'}$ and show how our $T$-deformations relate to the toric deformations constructed by K. Altmann in \cite{MR1329519} and \cite{MR1798979}.

We first briefly recall this construction. Consider some $k\in \NN$ and $R\in M'\setminus\{0\}$ and let $N=N'\cap R^\perp$ with cosection $s\colon N'\to N$, as in remark \ref{rem:downgrade}. We take $Q(R):=s(\sigma\cap [R=1])$ and $Q(R)^\infty:=s(\sigma\cap[R=0])$. Consider some admissible Minkowski decomposition
$Q(R)=Q_0+k\cdot Q_1+\ldots + k\cdot Q_r$
where if $k\neq 1$, $Q_1,\ldots, Q_r$ are lattice polyhedra.

For $R\in\sigma^\vee$, Altmann constructs a deformation of $X$ as follows. Take $\widetilde N=N\oplus \ZZ^{r+1}$, and consider 
the cone $\widetilde \sigma\subset \widetilde N_\QQ$ generated by 
$Q(R)^\infty\times 0$ and $Q_i\times e_i$ for $i=0,\ldots,r$, where $\{e_i\}$ is the standard basis of $\ZZ^{r+1}$.
We then get a map $\pi$ from $\tv(\widetilde{\sigma})$ to $\Aff^r$ by sending the coordinate function $t_i$, $1\leq i \leq r$, on $\Aff^r$ to $\chi^{k\cdot e_0^*}-\chi^{e_i^*}$, where $\{e_j^*\}$ is the dual basis of $\{e_j\}$, and for $u\in \widetilde M$, $\chi^{u}$ denotes the corresponding character.
It turns out that $\pi$ is flat with $\pi^{-1}(0)$ canonically isomporphic to $\tv(\sigma)$.

For $R\notin\sigma^\vee$, the construction is more complicated since the total space isn't toric. Consider $\tau:=
\sigma\cap [R\geq 0]$. Then there is a natural map $\tv(\tau)\to \tv(\sigma)$ presenting $\tv(\tau)$ as an open subset of a proper birational modification of $\tv(\sigma)$. Since $R\in\tau^\vee$, we can apply the above construction to get a deformation of $\tv(\tau)$. This deformation then maps birationally to Altmann's toric deformation of $\tv(\sigma)$, see \cite{MR1798979} theorem 3.2 for more details. 

\begin{remark}
For fixed $k$ and $R$, a toric deformation as above is called \emph{homogeneous of degree $-kR$}. This is justified by the fact that the image of such a deformation in $T_X^1$ under the Kodaira-Spencer map lies in the degree $-kR$ compoment of $T_X^1$ with respect to the natural $M'$-grading.  Furthermore, the images of toric deformations completely span $T_X^1$.
\end{remark}

To compare  toric deformations with our deformations, we view $X$ as a $T$-variety for the subtorus
$T = \ker (R \colon T' \to \CC^*) \subset T'$. As mentioned in remark \ref{rem:downgrade}, the corresponding 
polyhedral divisor $\D$ has coefficients
$$
\D_0 = s(\sigma \cap [R=1]) \qquad \D_\infty = s(\sigma \cap [R = -1])
$$
on $Y = \PP^1$, where $s\colon N'\to N$ is a cosection. Thus, $\D_0$ is Altmann's $Q(R)$ from above, and  the tail cone is
$Q(R)^\infty$.

In the case $R \in \sigma^\vee$, the coefficient at $\infty$ is empty, and
for admissible $r$-parameter Minkowski decompositions of $\D_0$, we
recover Altmann's homogeneous toric deformations in degree $-R$.
If $y = y_0$ is the coordinate on $\loc \D = \Aff^1$, a linear change
of coordinates $t_{0,i} \mapsto t_{0,i} + y$ moves the supporting prime
divisors of $\D^\tot$ on $\loc \D^\tot$ into the invariant divisors of the toric variety $\Aff^{r+1}$
so we see that $X^\tot$ is toric.
We can also see deformations in non-primitive degrees by considering decompositions of the form
$Q(R) = Q_0 + k\cdot Q_1+ \dotsc + k\cdot Q_r$ for some $k>1$  and applying remark \ref{rem:nonprim}, where we now require that $Q_i$ is a lattice polyhedron for $1\leq i \leq r$. 

In the case $R \not\in \sigma^\vee$, the total space of a $T$-deformation cannot be toric. Indeed, if $X^\tot$ were toric, by the downgrading procedure described in remark \ref{rem:downgrade}, $\loc \D^\tot$ would have to be a toric variety, with the prime divisors in the support of $\D^\tot$ all being invariant. As an example, consider the case $r = 1$ illustrated in
figure~\ref{fig:p1a1}. While $\PP^1 \times \Aff^1$ is toric, the union
of the three toric divisors $\{0\} \times \Aff^1$, $\{\infty\} \times \Aff^1$
and $\PP^1 \times \{0\}$ is connected. Since this is not the case for the three
divisors $D(0,0)$, $D(0,1)$ and $D(\infty,0)$, $X^\tot$ cannot have the structure of a toric variety. Here, we recover the
families which Altmann constructs by
deforming birational modifications of $X$.

\begin{remark}
While we saw above that $X^\tot$ is not necessarily toric,
it does retain an
action by  the larger torus $T'$. We will show this for the case of figure~\ref{fig:p1a1}; the general cases are similar.
Consider the diagonal action of $\CC^*$ on $\Yt = \PP^1 \times \Aff^1$. 
The prime divisors in $Y^\tot$ on which $\D^\tot$ is supported are all $\CC^*$-invariant. Thus, the $\CC^*$ acting on $Y^\tot$ also acts on $X^\tot$. Furthermore, $Y\subset Y^\tot$ is invariant under this $\CC^*$ action. Thus, $T'\cong T\times \CC^*$, and the embedding of $X$ in $X^\tot$ is equivariant. 

We can even describe  $X^\tot$ as a $T'$-variety in terms of a proper polyhedral divisor $\E^\tot$ on $\PP^1$. Indeed, we take
$$
\E^\tot=\E^\tot_0\otimes\{0\}+\E^\tot_1\otimes\{1\}+\E^\tot_\infty\otimes\{\infty\}
$$
with
\begin{align*}
\E^\tot_0 &= \alpha(\D^0_0 \times \{1\}) + \sigma' \\
\E^\tot_1 &= \alpha(\D^1_0 \times \{0\}) + \sigma' \\
\E^\tot_\infty &= \alpha(\conv\left\{\D_\infty \times \{-1\},
		\delta\times\{0\}\right\}) + \sigma',
\end{align*}
where $\alpha\colon N \oplus \ZZ \to N'$ is the isomorphism induced by $s$, and
$\sigma' = \sigma \cap [R \ge 0]$ is the positive part of $\sigma$.
We leave it to the reader to check that this gives us the total space $X^\tot=X(\E^\tot)$ as a $T'$-variety. Note that we still recover the Minkowski decomposition at this level:
$$
\E^\tot_0 + \E^\tot_1 = \sigma \cap [R \ge 1] \qquad
\E^\tot_\infty = \sigma \cap [R \ge -1].
$$
\end{remark}

\begin{figure}[htbp]
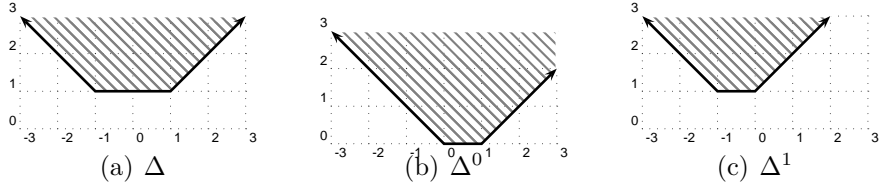

    \centering
    \subfigure[$\Delta$]{\fanopoly}
    \subfigure[$\Delta^0$]{\fanopolyb}
    \subfigure[$\Delta^1$]{\fanopolya}
    \caption{Minkowski decomposition for an affine threefold singularity}\label{fig:affineex}
\end{figure}

\begin{ex}
	We consider an example of a toric threefold with deformations in non-negative degrees. Let $N'=\ZZ^3$ with standard basis $e_1,e_2,e_3$ and $\sigma$ generated by $(-1,1,1)$, $(1,1,1)$, $(-1,1,-1)$, and $(1,1,-1)$. $X=\tv(\sigma)$ is then the cone over the singular projective Fano surface $X'$ whose minimal resolution is the toric surface presented in figure~\ref{fig:blowup} in section~\ref{sec:tvar}. 

	Setting $N=\langle e_1,e_2\rangle$ with cosection $s\colon N'\to N$ given by projection, we can consider $X$ as the $T$-variety $X(\D)$ over $Y=\mathbb{P}^1$ with $\D=\Delta\otimes \{0\}+\Delta\otimes\{\infty\}$ and $\Delta$ as in figure~\ref{fig:affineex}(a). The Minkowski decompositions $\D_0=\Delta^0+\Delta^1$ and $\D_\infty=\Delta^0+\Delta^1$ induce a two-parameter deformation $\pi$ of $X$. Restricting to the coordinate axes of the base space gives homogeneous deformations in degrees $-e_3^*$ and $e_3^*$, neither of which lie in $\sigma^\vee$.

	Note that the deformation $\pi$ has degree zero with respect to the $\ZZ=\langle e_2^* \rangle $ grading on $\CO_X$ inducing the quotient $X'$. Thus, $\pi$ induces a two-parameter deformation $\pi'$ on $X'$ as well.
\end{ex}

\section{Decompositions of Divisorial Fans}\label{sec:fandecomp}
Let $Y=\mathbb{P}^1$ and let $\dfan$ be a divisorial fan on $Y$. We now show how to construct deformations of the rational non-affine $T$-variety $X(\dfan)$. For simplicity's sake, we will restrict to those deformations which correspond to primitive degrees in the toric case.

\begin{defn}\label{def:mdcomp}
	Let $C$ be any (nonempty) polyhedral subdivision in $N_\QQ$ and $r\in\NN$. An admissible \emph{$r$-term Minkowski decomposition} of $C$ consists of admissible $r$-term Minkowski decompositions
$$
\Delta=\Delta^0+\Delta^1+\ldots+\Delta^r
$$
for all $\Delta\in C$ such that	
\begin{enumerate}
			\item \label{item:compatible} If $\Delta\cap\nabla\neq \emptyset$ with $\Delta,\nabla\in C$, then $(\Delta\cap\nabla)^{i}=\Delta^i\cap\nabla^i$ for any $i\in\{0,\ldots,r\}$.
	\item \label{item:complex}We have
				$$\sum_{i\in I}\bigcap_{\Delta\in\I} \Delta^i\prec\sum_{i\in I}\bigcap_{\Delta\in\J} \Delta^i  $$
				for any
				$\J\subset\I\subset C$
				and $I\subset\{0,\ldots,r\}$.\label{item:compface}
		\end{enumerate}
	\end{defn}

\begin{remark}
The above definition guarantees that for fixed $i$, the set $$C^i:=\{\Delta^i\ |\ \Delta\in C\}$$ is a polyhedral subdivision in $N_\QQ$. The admissibility of the decompositions of all $\Delta$ in $C$ is equivalent to the condition that for any vertex $v$ of $C$, at most one of the corresponding vertices $v^i\in C^i$ is not a lattice point.
As a slight abuse of notation, we will write $C=C^0+\ldots+C^r$ for a Minkowski decomposition; note that the decomposition also includes the data of which polyhedra in the $C^i$ are added together.
\end{remark}

Similar to section~\ref{sec:affinecase}, let $\mathcal{P}\subset Y$ be a finite set of points in $Y$,  this time including all those points $P$ with $\DP$ nontrivial for some $\D\in\dfan$. 
For each $P\in\mathcal{P}\subset Y$,  consider an admissible $r_P$-parameter Minkowski decomposition of the slice $\dfan_{P}$ for some $r_P\in\ZZ_{\geq 0}$. We call this a \emph{decomposition} of the divisorial fan $\dfan$. 

As before set $r=\sum r_P$. For each $\D\in\dfan$, we get a decomposition of $\D$ by decomposing $\DP=\sum_{i=0}^{r_P} \DP^i$, where if $\DP\neq \emptyset$, $\DP^i$ comes from the decomposition of $\dfan_P$, and for $\DP= \emptyset$, we set $\DP^i= \emptyset$. From the construction in section~\ref{sec:affinecase}, for each $\D\in\dfan$ we thus get a proper polyhedral divisor $\D^\tot$ on $Y^\tot=\PP^1\times\base$, along with a deformation $\pi_\D:X(\D^\tot)\to\base$ of $X(\D)$. We wish to glue these deformations together to get a deformation of $X(\dfan)$.

 For any $\I\subset \dfan$ and any $\lambda\in\base$, we set 
\begin{align*}
	\D^{\I}=\bigcap_{\D\in\I} \D;\qquad
	\D^{\I,\tot}=\bigcap_{\D\in\I} \D^\tot;\qquad
	\D^{\I,(\lambda)}=\bigcap_{\D\in\I} \D^{(\lambda)}.\\
\end{align*}
Note that $\D^\I\in\dfan$.
We then set 
$$
\dfan^\tot=\left \{\D^{\I,\tot}\right\}_{\I\subset \dfan};
\qquad\dfan^{(\lambda)}=\left\{\D^{\I,(\lambda)}\right\}_{\I\subset \dfan}. 
$$

\begin{prop}\label{prop:isfan}
$\dfan^\tot$ is a divisorial fan on $Y^\tot$. Likewise, each $\dfan^{(\lambda)}$ is a divisorial fan on $Y=\PP^1$.
\end{prop}	

We will prove this proposition shortly, but first we wish to construct the corresponding deformation. As before, there is a natural regular map $\pi\colon X({\dfan^\tot})\to \base$ by lemma~\ref{lemma:regular}.

\begin{thm}\label{thm:tdef}
	The map $\pi\colon X(\dfan^\tot)\to \base$ gives a flat family with $\pi^{-1}(\lambda)\cong X(\dfan^{(\lambda)})$ for $\lambda\in \base$. In particular, $\pi^{-1}(0)=X(\dfan)$.
\end{thm}	
We call deformations of the type constructed in the above theorem \emph{$T$-deformations}.

\begin{figure}[htbp]
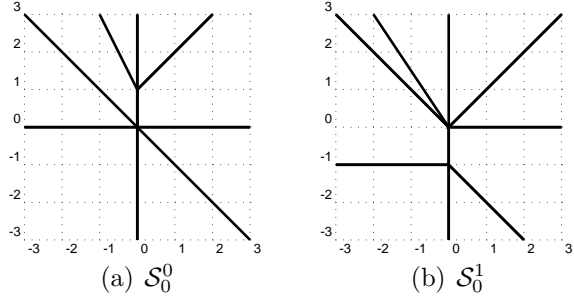

    \centering
    \subfigure[$\dfan_0^0$]{\cotangfannulla}
    \subfigure[$\dfan_0^1$]{\cotangfannullb}
    \caption{A deformation of $\mathbb{P}(\Omega_{\mathcal{F}_1})$}\label{fig:cotanf1def}
\end{figure}
\begin{ex}
	Consider the $T$-variety $X=X(\dfan)=\mathbb{P}(\Omega_{\mathcal{F}_1})$ as described in the example at the end of section~\ref{sec:tvar} and in figure~\ref{fig:cotanf1}. We construct a one-parameter deformation over $\mathbb{A}^1\setminus\{1\}$ by decomposing the slice $\dfan_0=\dfan_0^0+\dfan_0^1$ as pictured in figure~\ref{fig:cotanf1def}. The general fiber of the deformation is 
$X(\dfan^{(\lambda)})$ for general $\lambda$, where $\dfan_0^{(\lambda)}=\dfan_0^0$, $\dfan_\lambda^{(\lambda)}=\dfan_0^1$, $\dfan_1^{(\lambda)}=\dfan_1$, and  $\dfan_\infty^{(\lambda)}=\dfan_\infty$, such that any  two-dimensional polyhedra with common tailcones are coefficients for the same polyhedral divisor.

\begin{figure}[htbp]
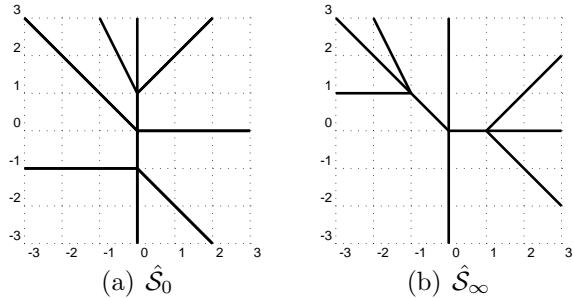

    \centering
    \subfigure[$\hat{\dfan}_0$]{\cotangfannull}
    \subfigure[$\hat{\dfan}_\infty$]{\cotangfaninftydegen}
    \caption{A degeneration of $\mathbb{P}(\Omega_{\mathcal{F}_1})$}\label{fig:cotanf1degen}
\end{figure}

On the other hand, $\mathbb{P}(\Omega_{\mathcal{F}_1})$ is the general fiber in a $T$-deformation of a toric variety. Indeed, if $\hat{\dfan}$ is the divisorial fan on $\mathbb{P}^1$ pictured in figure~\ref{fig:cotanf1degen} (where once again two-dimensional polyhedra with common tailcones belong to the same polyhedral divisor), $X(\hat{\dfan})$ is in fact a toric variety and $X(\dfan)$ is the general fiber of the homogeneous deformation coming from the decomposition $\hat{\dfan}_\infty=\dfan_\infty+\dfan_1$. This is similar to the degeneration of $\mathbb{P}(\Omega_{\mathbb{P}^2})$ to the projective cone over the del Pezzo surface of degree six constructed in~\cite{suess08}, example~5.1. 
 \end{ex} 

The remainder of the section is dedicated to proving proposition \ref{prop:isfan} and theorem \ref{thm:tdef}. We split up the proofs into several smaller lemmata. Note that we will only be proving the claim of proposition \ref{prop:isfan} for $\dfan^\tot$. The proof for $\dfan^{(\lambda)}$ is similar, and left to the reader.

\begin{lemma}\label{lemma:Qface}
	Consider polytopes $\nabla^i\subset \Delta^i$ for $1\leq i\leq n$. Set $\Delta:=\sum_{i=1}^n \Delta^i$ and $\nabla:=\sum_{i=1}^n \nabla^i$ and let $I$ be any subset of $\{1,\ldots,n\}$.
	\begin{enumerate}
		\item 	For any $w\in(\tail \Delta)^\vee$ with $\face( \Delta,w)=\nabla$, we have $\face(\sum_{i\in I} \Delta^i,w)=\sum_{i\in I} \nabla^i.$\label{item:Qface1}\\
		\item 	For any $w\in(\tail \Delta)^\vee$ with $\face( \Delta,w)=\face(\nabla,w)$, we have $\face(\sum_{i\in I} \Delta^i,w)=\face(\sum_{i\in I} \nabla^i,w).$\label{item:Qface2}\\
\end{enumerate}

\end{lemma}
\begin{proof}
Note that the first claim follows from the second. Indeed, if $\nabla=\face(\Delta,w)$, then $w$ is constant on $\nabla$, and must also be constant on $\nabla^i$. Thus, $\nabla^i=\face(\nabla^i,w)$.

For part \ref{item:Qface2}, observe that $\nabla^i\subset \Delta^i$ implies
$\min \langle \nabla^i,w\rangle \geq \min \langle \Delta^i,w\rangle$ for $1\leq i \leq n$. But 
$$\sum \min \langle \nabla^i,w\rangle=\min \langle \nabla,w\rangle=\min \langle \Delta,w\rangle=\sum \min \langle \Delta^i,w\rangle$$
so we in fact have $\min \langle \nabla^i,w\rangle = \min \langle \Delta^i,w\rangle$. Coupled with $\nabla^i\subset \Delta^i$ we then get that $\face(\nabla^i,w)\subset \face(\Delta^i,w)$.
Since taking faces commutes with Minkowski sums, we then have the following diagram:
$$
\begin{array}{c c c c c c c c c c c c c}
	\face(\Delta^{1},w)&+&\ldots&+&\face(\Delta^{n},w)&=& \face(\Delta,w)\\
\cup&&&&\cup&& \Vert\\
\face(\nabla^{1},w)&+&\ldots&+&\face(\nabla^{n},w)&=& \face(\nabla,w).\\
\end{array}
$$
We can conclude that the inclusions must be equalities. The claim then follows by again using the fact that taking faces commutes with Minkowski sums.
\end{proof}

\begin{lemma}\label{lemma:alwayspos}
Let $\D',\D$ be proper polyhedral divisors on some curve $C$ with $\D'\prec \D$, $\loc (\D)$ complete, and $\loc (\D')$ not complete. Then for any $w\in \tail (\D)^\vee$ with $$\face(\tail(\D),w)=\tail(\D'),$$ we have $\deg(\D(w))>0$.
\end{lemma}
\begin{proof}
	Let $y$ be the general point of $C$. Then since $\D'\prec \D$, we find $w_y$ such that $\face(\tail(\D),w_y)=\tail(\D')$ and $\deg(\D(w_y))>0$, since the fact that some coefficients of $\D'$ are the empty set implies the existence of a divisor $D_y\in|\D(w_y)|$ with nontrivial support. From this it follows that $\deg (\D)\cap \tail(\D')=\emptyset$.

	Now let $w$ be as in the statement of the lemma. The hyperplane determined by $\langle \cdot , w \rangle =0$ intersects $\tail(\D)$ in exactly $\tail(\D')$. Thus, it cannot intersect $\deg (\D)$, since $\deg(\D)\subset \tail(\D)$. It follows that $\deg (\D(w))=(\deg (\D))(w)\neq 0$, so it must be strictly positive.
\end{proof}

\begin{lemma}\label{lemma:iidef}
Consider any $\I\subset \dfan$.
\begin{enumerate}
\item  $\D^{\I,\tot}$ arises from an admissible Minkowski decomposition of $\D^\I$.\label{item:marise}
\item $\D^{\I,\tot}$ is proper.\label{item:iprop}
\end{enumerate}
\end{lemma}
\begin{proof}
Consider $P\in\mcP$, and set $\DP^{\I,i}:=\D_x^{\I,\tot}$ for $x=V(y_P-t_{P,i})$. Then clearly $\DP^\I=\sum \DP^{\I,i}$, and all $\DP^{\I,i}$ have the correct tail cone (or are the empty set). Thus, we just need to check the admissibility of the decomposition. But from definition \ref{def:mdcomp}\ref{item:complex} coupled with \ref{lemma:Qface}\ref{item:Qface1} we have that for any $\D\in\I$ there exists $w\in\tail(\D)^\vee$ such that for $0\leq i \leq r_P$ either $\DP^{\I,i}=\emptyset$ or $\DP^{\I,i}=\face(\DP^i,w)$. Thus, since the $\DP^i$ form an admissible decomposition, the $\DP^{\I,i}$ must as well, proving the first claim.

The second claim follows from the first coupled with lemma \ref{lemma:properness}. 
\end{proof}

The next lemma is the essential point in the proof of proposition \ref{prop:isfan}. It is rather technical in that we must consider a number of different cases, but each case just requires an application of some of the above lemmata.

\begin{lemma}\label{lemma:totface}
	For any $\J\subset\I\subset \dfan$, $\D^{\I,\tot}\prec \D^{\J,\tot}$.
\end{lemma}
\begin{proof}
	For any point $x\in Y^\tot$ not necessarily closed, we define a point $\hat{x}\in Y$ as follows. If $x$ is contained in some $V(y_P-t_{P,i})$ for $P\in\mcP$ and $0\leq i \leq r_P$, then let $\hat{x}$ be equal to the point $P$ in $\PP^1$. Note that such a $P$, if it exists, is unique due to the way $\base$ was constructed. In this case, we say $x$ is \emph{special}. Otherwise, if $x$ is contained in any divisor $V(y_P)$ for $P\in Y\setminus \mcP$, let $\hat{x}=P$. Finally, for any other point $x$, take $\hat{x}$ to be the general point in $\PP^1$. We now describe $\D_x^{\I,\tot}$. If $x$ is contained in some $V(y_P-t_{P,i})$, let $I$ be the set of all $i$ for which this holds. Then
	one easily sees that $\D_x^{\I,\tot}=\sum_{i\in I} \bigcap_{\D\in\I}\DP^i$. Otherwise, $\D_x^{\I,\tot}=\D_{\hat{x}}^\I$. One can describe $\D_x^{\J,\tot}$ similarly.

	Fix now any point $y\in \loc(\D^{\I,\tot})$, not necessarily closed. We must show that the requirements of definition \ref{def:face}\ref{item:face} hold for this point $y$. Now, since $\D^\I\prec \D^\J$, there exists $(w_{\hat{y}},D_{\hat{y}})\in M\times|\D^\J(w_{\hat{y}})|$ fulfilling the face relation of definition \ref{def:face} for $\D^\I\prec \D^\J$ and the point $\hat{y}$. In the remainder of the proof, we will consider several cases:
\begin{enumerate}[label=(\alph*)]
	\item $\loc (\D^{\I,\tot})=Y^\tot$ and $\deg (D_{\hat{y}})=0$;\label{item:tfc1}
	\item $\loc (\D^{\I,\tot})=Y^\tot$ and $\deg (D_{\hat{y}})>0$;\label{item:tfc2}
	\item $\loc (\D^{\I,\tot})\neq Y^\tot$ and $\loc (\D^{\J,\tot})= Y^\tot$;\label{item:tfc3}
	\item $\loc (\D^{\J,\tot})\neq Y^\tot$.\label{item:tfc4}
\end{enumerate}

Starting with case \ref{item:tfc1}, set $w_y=w_{\hat{y}}$. We take $D_y$ to be the trivial divisor on $Y^\tot$. Note that we have $D_y\in|\D^{\J,\tot}(w_y)|$. Clearly $y\notin \supp D_y$. Furthermore, we claim 
\begin{equation}
	\face(\D_y^{\J,\tot},w_y)=\D_y^{\I,\tot}.\label{eqn:face1}
\end{equation}
Indeed, if $y$ isn't special, then this follows from $\D_y^{\I,\tot}=\D_{\hat{y}}^\I$ and $\D_y^{\J,\tot}=\D_{\hat{y}}^\J$. 
For $y$ special, point \ref{item:compatible} of definition \ref{def:mdcomp} gives us
\begin{align*}
	\DQ^\I=\sum_{i=0}^{r_Q}  \bigcap_{\D\in\I}(\DQ^i)\\
	\DQ^\J=\sum_{i=0}^{r_Q} \bigcap_{\D\in\J}(\DQ^i)
\end{align*}
whereas we automatically have
$$
\bigcap_{\D\in\I}\DQ^i\subset \bigcap_{\D\in\J}\DQ^i
$$
for all $0\leq i \leq r_Q$.
Since $\face(\DQ^\J,w_y)=\DQ^\I$, we can thus apply \ref{lemma:Qface}\ref{item:Qface1} to show equation \eqref{eqn:face1}.
Now finally, for all $v\in Y^\tot$, we claim that  
$
\face(\D_v^{\J,\tot},w_y)=\face(\D_v^{\I,\tot},w_y).
$
Indeed, for all $v$ we have $\face(\D_v^{\J},w_y)=\face(\D_v^{\I},w_y)$. For $v$ not special the claim is then immediate. On the other hand, for $v$ special we use lemma 
\ref{lemma:Qface}\ref{item:Qface2}, where the hypothesis of the lemma is once again satisfied due to point \ref{item:compatible} of definition \ref{def:mdcomp}. Thus, the pair $(w_y,D_y)$ satisfies the requirements of definition \ref{def:face}. 

We now move to case \ref{item:tfc2}. Since $\deg (D_{\hat{y}})>0$, clearly we can find some $k\in\NN$ such that $|\D^{\J,\tot}(k\cdot w_{\hat{y}})|$ contains a divisor $D$ such that $Y^\tot\setminus\supp D$ contains only those special points $x$ with $\hat{x}=\hat{y}$. We then set $w_y=k\cdot w_{\hat{y}}$ and take $D_y=D$. Now, we have $y\notin \supp D_y$, and
$
\face(\D_y^{\J,\tot},w_y)=\D_y^{\I,\tot}
$
exactly as in case \ref{item:tfc1}. We claim that we also have
$
\face(\D_v^{\J,\tot},w_y)=\face(\D_v^{\I,\tot},w_y)
$
for all $v\in Y^\tot\setminus\supp D_y$. If $\hat{v}\neq \hat{y}$, then  $\D_v^{\I,\tot}$ and $\D_v^{\J,\tot}$ are both trivial and the claim follows from the properties of $w_{\hat{y}}$. Likewise, if $\hat{v}= \hat{y}$ and $y$ is trivial, then $\D_v^{\I,\tot}=\D_y^{\I}$, $\D_v^{\J,\tot}=\D_y^{\J}$ and the claim again follows from the properties of $w_{\hat{y}}$. Finally, if $\hat{v}= \hat{y}$ and $y$ is not trivial, we can apply lemma \ref{lemma:Qface}\ref{item:Qface2} as in part \ref{item:tfc1}. Thus, we again have that the pair $(w_y,D_y)$ satisfies the requirements of definition \ref{def:face}. 

We now consider case \ref{item:tfc3}. Suppose first that $\hat{y}\in \loc (\D^\I)$. Then one easily sees that $\deg (D_{\hat{y}})>0$ and one can proceed as in case \ref{item:tfc2}. Thus, we have reduced to the case that $\hat{y}\notin \loc( \D^\I)$, from which it follows that $y$ must be special. Now, we can find $w\in\tail(\D^\J)$ with $\face(\D_y^{\J,\tot},w)=\D_y^{\I,\tot}$ by point 
\ref{item:complex} of definition \ref{def:mdcomp}. Furthermore, by lemma \ref{lemma:alwayspos} we have $\deg (\D^{\J}(w))>0$ since taking faces and tails commutes.
Similar to in case \ref{item:tfc2},
we can find some $k\in\NN$ such that $|\D^{\J,\tot}(k\cdot w)|$ contains a divisor $D$ such that $Y^\tot\setminus\supp D$ contains only those special points $x$ with $\hat{x}=\hat{y}$ and $\D_x^{\I,\tot}\neq \emptyset$. We then set $w_y=k\cdot w_{\hat{y}}$ and take $D_y=D$.
The claim of 
$
\face(\D_y^{\J,\tot},w_y)=\D_y^{\I,\tot}
$
is satisfied automatically by our choice of $w_y$.
The claim that 
$
\face(\D_v^{\J,\tot},w_y)=\face(\D_v^{\I,\tot},w_y)
$
for all $v\in Y^\tot\setminus\supp D_y$ is immediate for $v$ nonspecial, and follows from lemma \ref{lemma:Qface}\ref{item:Qface2} for $v$ special. Thus, we again have that the pair $(w_y,D_y)$ satisfies the requirements of definition \ref{def:face}. 

The final case \ref{item:tfc4} is essentially identical to the case of \ref{item:tfc3}, with the simplification that since $\loc (\D^\J)$ is affine, we needn't worry about the degrees of evaluations of $\D^\J$.
\end{proof}

We are now ready to prove proposition \ref{prop:isfan}:
\begin{proof}[Proof of proposition \ref{prop:isfan}]
	First, all elements of $\dfan^\tot$ are proper polyhedral divisors due to lemma \ref{lemma:iidef}\ref{item:iprop}. Secondly, we claim that intersections of elements of $\dfan^\tot$ are themselves elements of $\dfan^\tot$. Indeed, for $\I,\J\subset\dfan$, we have $\D^{\I,\tot}\cap\D^{\J,\tot}=\D^{\I\cup\J,\tot}$. Finally, from lemma \ref{lemma:totface} we have the necessary face relations:
	$$
	\D^{\I,\tot}\succ\D^{\I\cup\J,\tot}\prec \D^{\J,\tot}
	$$
\end{proof}

We conclude the section with the proof of theorem \ref{thm:tdef}:
\begin{proof}[Proof of theorem \ref{thm:tdef}]
	Since the maps $\pi_\D$ arise from a projection of the quotient map, they  agree along intersections of polyhedral divisors and we can clearly glue them together to the map $\pi$.
	Flatness of $\pi$ can be checked locally on each $X(\D^{\I,\tot})$ for $\I\subset\dfan$; this follows then directly from lemma \ref{lemma:iidef}\ref{item:marise} and theorem \ref{thm:affinemainthm}. From this theorem, we also know $\pi_{|X(\D^{\I,\tot})}^{-1}(\lambda)=X(\D^{\I,(\lambda)})$, so we just need to check that everything glues properly. But for $\I,\J\subset \dfan$, 
	$$\pi_{|X(\D^{\I\cup\J,\tot})}^{-1}(\lambda)=X(\D^{\I\cup\J,(\lambda)})=X(\D^{\I,(\lambda)})\cap X(\D^{\J,(\lambda)})=\pi_{|X(\D^{\I,\tot})}^{-1}(\lambda)\cap\pi_{|X(\D^{\J,\tot})}^{-1}(\lambda).$$
	Thus, the gluing on $X(\dfan^\tot)$   induces the gluing on $X(\dfan^{(\lambda)})$.
	\end{proof}

\section{Locally Trivial Deformations}\label{sec:ks}
Let $Y=\mathbb{P}^1$ and let $\dfan$ be a divisorial fan on $Y$. To any locally trivial one-parameter deformation of $X(\dfan)$ we can assign a class in $H^1(X(\dfan),\T_{X(\dfan)})$ via the Kodaira-Spencer map: we pull back the deformation to one over $\spec \CC[t]/t^2$, and the Kodaira-Spender correspondence gives a bijection between such first-order deformations and the above-mentioned cohomology classes. In this section, we will compute the image of this map  for certain special $T$-deformations.

Now for some $P\in Y$ let $\dfan_P=\dfan_P^0+\dfan_P^1$ be an admissible one-parameter Minkowski decomposition of $\dfan_P$ such that for each $\Delta\in\dfan_P$, either $\Delta^0$ or $\Delta^1$ is a lattice translate of $\tail(\Delta)$. Note that this is always the case if $X(\dfan)$ is smooth; indeed this follows from \cite{liendo:10a} proposition~5.1 and theorem~5.3 together with the fact that height-one hyperplane sections of smooth cones only have trivial admissible Minkowski decompositions. We furthermore assume that for any $\D\in\dfan$ with $\loc(\D)=Y$, $\D$ has trivial coefficients everywhere except $P$ and $\infty:=V(y_P^{-1})$. The calculations can be carried out under more general assumptions, but the resulting formula is considerably more complicated, see \cite{ilten:10a} theorem 3.5.1.

We consider the one-parameter deformation $\pi\colon X(\dfan^\tot)\to \base$ corresponding to the above decomposition. We shall describe its image in $H^1(X(\dfan),\T_{X(\dfan)})$ as a \v{C}ech cocycle and thus need to choose some open cover $\mathfrak{U}$ of $X(\dfan)$. By possibly refining the original divisorial fan $\dfan$, we can assume without loss of generality that for $\D\in \dfan$ with affine locus, either $\DP=\emptyset$, or $\DQ=\tail(\DQ)$ for all $Q\in\loc \D$, $Q\neq P$. We will then use the open cover $\mathfrak{U}=\{X(\D)\ |\ \D\in\dfan\}$ provided by this divisorial fan to describe \v{C}ech cocycles.

To the above Minkowski decomposition and for every $\D\in \dfan$, we can associate $a_\D\in\{1,-1\}$ and $\lambda_\D\in N$ as follows. Consider first $\D\in\dfan$ with $\DP\neq \emptyset$.  
If we can write $\DP=\lambda_\D+{\DP^{0}}$ for some $\lambda_\D\in N$, let $a_\D=1$. Otherwise set $a_\D=-1$ and define $\lambda_\D$ by $\DP=\lambda_\D+{\DP^{1}}$. Note that one of these conditions must be fulfilled due to our above assumption on the nature of the decomposition of $\dfan_P$.  
For $\D\in\dfan$ with $\DP=\emptyset$, we take $a_\D=1$ and $\lambda_\D=0$. Finally, let $e_i^*$ be a basis for $M$ and take $t=t_{P,1}$.

 \begin{thm}\label{thm:tvarks}
	The deformation $\pi\colon X(\dfan^\tot)\to S$ is locally trivial after pullback to $\spec \CC[t]/t^2$ and its image under the Kodaira-Spencer map is the cocycle defined by
$$d_{\D,\E}=\frac{a_\D-a_\E}{2}\frac{\partial}{\partial y_P}+y_P^{-1}\sum_i \langle a_\D\lambda_\D-a_\E\lambda_\E,e_i^*\rangle \chi^{e_i^*} \frac{\partial}{\partial \chi^{e_i^*}} $$
for $\D,\E\in\dfan$.
\end{thm}

For the proof of the theorem, we shall use the following two lemmata, which are essentially special cases of proposition~8.6 in~\cite{MR2207875}:
\begin{lemma}\label{lemma:latticetrafo}
	Let $Y$ be a smooth variety and $\dfan$ a divisorial fan on $Y$. For some $v\in N$ and $D_0,D_1\in\WDiv(Y)$ with $D_1-D_0=\Div(f)$ let $\widetilde{\dfan}=\{\D+v(D_1-D_0)\mid\D\in\dfan\}$. Then there is a canonical isomorphism $\phi_v\colon X(\dfan)\rightarrow X(\widetilde{\dfan})$ where $\phi_v^\#$ is defined by mapping $\chi^u$ to $f^{\langle v,u\rangle}\chi^u$ for $u\in M$.
	\begin{proof}
	We show that locally $\phi_v$ is an isomorphism. Consider $\D \in \dfan$ with tailcone $\delta$ and set $\widetilde{\D}=\D+v(D_1-D_0)$. Then 
	\begin{align*}
		\CO_{X(\widetilde{D})}=\bigoplus_{u\in \delta^\vee\cap M} H^0\Big(Y,\widetilde{\D}(u)\Big)\cdot \chi^u
		&\cong \bigoplus_{u\in \delta^\vee\cap M} H^0\Big(Y,\D(u)+ \langle v,u\rangle\Div(f)\Big)\cdot \chi^u\\
		&\cong \bigoplus_{u\in \delta^\vee\cap M} H^0\Big(Y,\D(u)\Big)\cdot f^{-\langle v,u\rangle}\chi^u\\
	\end{align*}
and thus $\phi_v^\#$ induces an isomorphism $\CO_{X(\widetilde{D})}\cong\CO_{X({D})}$, since $$\CO_{X({D})}=\bigoplus_{u\in \delta^\vee\cap M} H^0\Big(Y,\D(u)\Big)\cdot \chi^u.$$
	\end{proof}
\end{lemma}

\begin{lemma}\label{lemma:baseauto}
	Let $Y$ be a smooth variety and $\bar{\gamma}\in\aut(Y)$. For a proper polyhedral divisor $\D$ on $Y$, define $\bar{\gamma}(\D)=\sum_D \D_D\cdot \bar{\gamma}(D)$.  Then there is a natural isomorphism $\gamma\colon X(\D)\to X(\bar{\gamma}(\D))$ induced by $\bar{\gamma}$.
	\begin{proof}
		Similar to the proof of the above lemma, we have
		\begin{align*}
			\CO_{X(\bar{\gamma}(\D))}=&\bigoplus_{u\in \delta^\vee\cap M} H^0\Big(Y,\bar{\gamma}({\D})(u)\Big)\cdot \chi^u\\
			=&\bigoplus_{u\in \delta^\vee\cap M} (\bar{\gamma}^\#)^{-1}\Big(H^0\big(Y,{\D}(u)\big)\Big)\cdot \chi^u.
		\end{align*}

	\end{proof}
\end{lemma}

To deal with polyhedral divisors with affine locus, we also need the following lemma:
	\begin{lemma}\label{lemma:emptyiso}
		Consider two proper polyhedral divisors $\E,\F$ on $Y^\tot$ supported on divisors of the form $V(y_Q-ct)$ for $Q\in Y$, $c\in\CC$. Suppose that for any prime divisor $D=V(y_Q-ct)\subset Y^\tot$ of the above form, either $\E_D=\F_D$, or $(\E_{|Y})_Q=(\F_{|Y})_Q=\emptyset$. Then $X(\E)\times_\base \spec\CC[t]/t^2=X(\F)\times_\base \spec\CC[t]/t^2$.
	\end{lemma}
	\begin{proof}
	 If  $\E$ or $\F$ has locus $Y^\tot$ then the statement is trivial. If not, then the loci are in fact affine. 
The claim follows from the fact that if $(\E_Y)_Q=\emptyset$, then
$$
(y_Q-ct)^{-1}\equiv y_Q^{-1}(1+cy_Q^{-1}t) \quad (t^2)$$
is regular on $\loc \E \times_\base \spec\CC[t]/t^2$ for all $c\in\CC$.
\end{proof}

\begin{proof}[Proof of theorem \ref{thm:tvarks}]
	We calculate the image of the Kodaira-Spencer map as described in~\cite{MR2247603}. First, we show that $\pi$ pulled back to $\spec \CC[t]/t^2$ is locally trivial with respect to the cover $\mathfrak{U}$ by constructing isomorphisms $$\begin{CD}\theta_\D \colon X(\D)\times\spec \mathbb{C}[t]/t^2@>\sim>> X(\D^\tot)\times_{\base}\spec \mathbb{C}[t]/t^2\end{CD}$$ 
for each $\D\in\dfan $.
	Let $\overline{\gamma}_{-1}\colon Y^\tot\to Y^\tot$ be the automorphism induced by $\bar{\gamma}_{-1}^{\#}\colon t\mapsto t$ and  $\bar{\gamma}_{-1}^{\#}\colon y_P\mapsto y_P+t$. Likewise, let $\overline{\gamma}_{1}\colon Y^\tot\to Y^\tot$ be the identity. Furthermore, let $\gamma_i$, $i=-1,1$ be the corresponding morphisms from lemma~\ref{lemma:baseauto}. We also let $\phi_v$ be as in lemma~\ref{lemma:latticetrafo} where we take divisors $D_0=V(y_P)$, $D_1=V(y_P-t)$ on $Y^\tot$. We claim that setting
	$$\theta_\D=\phi_{a_\D\lambda_\D}\gamma_{a_\D}$$
	gives the desired isomorphism after taking the fiber product with $\spec \CC[t]/t^2$ over $\base$.

To begin with, we have that
	$$X(\D)\times\base =X(\sum_Q\D_{Q}\otimes D_Q)
$$
maps via $\gamma_{a_\D}$ to
\begin{equation}\label{eqn:uglything}
X\left(\D_{P} \otimes \left(D_0+\frac{a_D-1}{2}(D_0-D_1)\right)+\D_\infty\otimes D_\infty+\sum_{Q\neq P,\infty} \DQ\otimes \bar\gamma_i(D_Q)\right)
\end{equation}
where $D_Q=V(y_Q)$.
This is equal to 
\begin{equation}\label{eqn:zs}
X\left(\D_{P} \otimes \left(D_0+\frac{a_\D-1}{2}(D_0-D_1)\right)+\D_\infty\otimes D_\infty +\sum_{Q\neq P,\infty} \DQ\otimes D_Q \right)
\end{equation}
after taking the fiber product with $\spec \CC[t]/t^2$ over $\base$.
Indeed, for $\loc \D$ complete or $a_\D=1$, equality holds outright. For $\loc \D$ affine and $a_\D=-1$, we simply apply lemma \ref{lemma:emptyiso}.

We then apply $\phi_{a_\D\lambda_\D}$ to \eqref{eqn:zs} to get 
$$
X\left(\DP^{0}\otimes D_0+\DP^{1}\otimes
 D_1+ \sum_{Q\neq P} \DQ\otimes D_Q\right)=X(\D^\tot).
$$
Thus, $\theta_\D$ gives the desired isomorphism after taking the fiber product with $\spec \CC[t]/t^2$ over $\base$.

Now set $\theta_{\D,\E}=\theta_\D^{-1} \theta_\E$.
	Define a derivation $d_{\D,\E}$ via $\theta_{\D,\E}^{\#}=\id{} + t \cdot d_{\D,\E}$. 
	For simplicity of notation, set $b_\D=\frac{a_\D-1}{2}$.	Now $\theta_{\D,\E}^{\#}=\gamma_{a_\E}^{\#}\phi_{a_\E\lambda_\E-a_\D\lambda_\D}^{\#}(\gamma_{a_\D}^{\#})^{-1}$ and we can calculate
\begin{align*}
	\gamma_{a_\E}^{\#}\phi_{a_\E\lambda_\E-a_\D\lambda_\D}^{\#}(\gamma_{a_\D}^{\#})^{-1}\left(y_P\right)&=\gamma_{a_\E}^{\#}\phi_{a_\E\lambda_\E-a_\D\lambda_\D}^{\#}\left(y_P+b_\D t\right)=\gamma_{a_\E}^{\#}\left(y_P+ b_\D t\right)\\&=y_P+(b_\D-b_\E)t=y_P+\frac{a_\D-a_\E}{2}\cdot t\\
	\end{align*}
		and
\begin{align*}
\gamma_{a_\E}^{\#}\phi_{a_\E\lambda_\E-a_\D\lambda_\D}^{\#}(\gamma_{a_\D}^{\#})^{-1}\left(\chi^{u}\right)
=\gamma_{a_\E}^{\#}\phi_{a_\E\lambda_\E-a_\D\lambda_\D}^{\#}\left(\chi^{u}\right)
=\gamma_{a_\E}^{\#}\left(\left(\frac{y_P-t}{y_P}\right)^{\langle a_\E\lambda_\E-a_\D\lambda_\D,u\rangle}\chi^u\right)\\
=\left(\gamma_{a_\E}^{\#}\left(\frac{y_P-t}{y_P}\right)\right)^{\langle a_\E\lambda_\E-a_\D\lambda_\D,u\rangle}\chi^u
=\left(\frac{y_P-(1+b_j)t}{y_P-b_jt}\right)^{\langle a_\E\lambda_\E-a_\D\lambda_\D,u\rangle}\chi^u\\
=\left(\frac{y_P-a_\E t}{y_P}\right)^{a_\E\langle a_\E\lambda_\E-a_\D\lambda_\D,u\rangle}\chi^u=\left(1-a_\E y_P^{-1}t\right)^{a_\E\langle a_\E\lambda_\E-a_\D\lambda_\D,u\rangle}\chi^u,
		\end{align*}
		where the equality at the start of the last line can be seen by considering the cases $a_\E=1$ and $a_\E=-1$ separately.
Thus,
\begin{align*}
d_{\D,\E}(y_P)&=\frac{a_\D-a_\E}{2}\\
d_{\D,\E}(\chi^u)&=\langle a_\D\lambda_\D-a_\E\lambda_\E,u\rangle \cdot y_P^{-1}\chi^u\end{align*}
and the theorem follows.
\end{proof}

\begin{remark}
Assume $X=X(\dfan)$ is smooth, and suppose that for some point $P\in\PP^1$, if for any $\D\in\dfan$ $\loc D$ is complete, then $\DQ$ is trivial for all $Q\neq P,\infty$, where $\infty=V(y_P^{-1})$.
	Although the above theorem is only stated for one-parameter deformations, it can be used to calculate the linear Kodaira-Spencer map $T_{\base,0}\to H^{1}(X,\T_X)$ for an $r$-parameter deformation coming from a decomposition $\dfan_P=\dfan_P^0+\ldots+\dfan_P^r$. Indeed, $T_{\base,0}=T_{\mathbb{A}^k,0}$ and each basis vector of the natural basis of $T_{\mathbb{A}^k,0}$ corresponds to a one-parameter deformation as described above.
\end{remark}

\section{Deformations of Non-affine Toric Varieties}\label{sec:tc}
As mentioned in the remark at the end of section~\ref{sec:tvar}, a toric variety can be viewed as a $T$-variety by considering the action of some subtorus of the big torus. Thus, we can use the divisorial fan decompositions of section~\ref{sec:fandecomp} to construct deformations of an arbitrary toric variety. If in particular the deformation is locally trivial, we can use theorem~\ref{thm:tvarks} to calculate the image of the Kodaira-Spencer map. We will use this result to show that  $T$-deformations span the space of first-order deformations for any complete smooth toric variety.

Let $N'$ be an $n$-dimensional lattice with dual $M'$; choose some basis $e_1,\ldots, e_n$ of $N'$ with corresponding dual basis $e_1^*,\ldots, e_n^*$. Let $\Sigma$ be a fan on $N_\QQ'$ with corresponding toric variety $X=\tv(\Sigma)$. We can consider $X$ as a $T$-variety with respect to the subtorus $T^N$, where $N$ is generated by the first $n-1$ basis elements and the cosection $s\colon N'\to N$ is simply the natural projection. In this setting, a $T$-deformation of $X$ is given by Minkowski decompositions of the polyhedral complexes $\dfan_0=s(\Sigma\cap [e_n^* =1])$ and $\dfan_\infty=s(\Sigma\cap[ e_n^* =-1])$. Here, we always take $y_0=\chi^{e_n^*}$ and $y_\infty=\chi^{-e_n^*}$.

\begin{figure}[htbp]
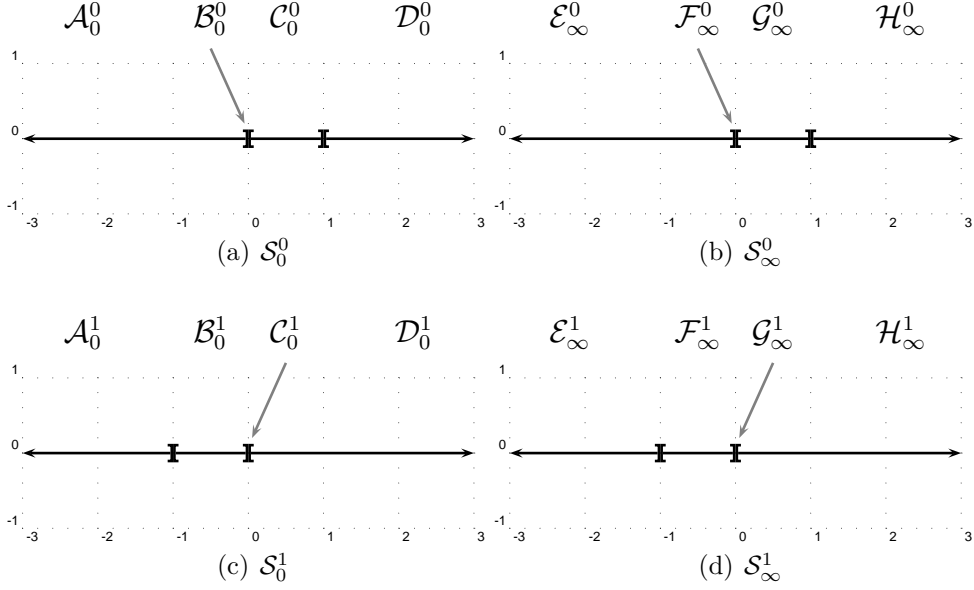

    \centering
    \subfigure[$\dfan_0^0$]{\blowupnulla}
    \subfigure[$\dfan_\infty^0$]{\blowupinftya}\\
    \subfigure[$\dfan_0^1$]{\blowupnullb}
    \subfigure[$\dfan_\infty^1$]{\blowupinftyb}\\

\caption{Slice decompositions for a toric surface}\label{fig:blowupdef}
\end{figure}

\begin{ex}
	Consider the toric surface $X$ in the first example in section~\ref{sec:tvar}, see figure \ref{fig:blowup}. A two parameter deformation of $X$ can be constructed by decomposing $\dfan_0=\dfan_0^0+\dfan_0^1$ and  $\dfan_\infty=\dfan_\infty^0+\dfan_\infty^1$ as pictured in figure \ref{fig:blowupdef}. Note that the coefficients $\B_0^0,\C_0^1, \F_\infty^0, \G_\infty^1$ are all simply the lattice point $0$.  In this case, the base space $\base\subset \mathbb{A}^2$ is the complement of the set $V(t_{0,1} t_{\infty,1}-1)$. 
\end{ex}

Now take some  one-parameter $T$-deformation $\pi$ of $X$ coming from a Minkowski decomposition  $\dfan_0=\dfan_0^0+\dfan_0^1$. Suppose additionally that the decomposition of $\dfan_0$ is such that for all $\D\in\dfan$ with $\D_0\neq \emptyset$, $\D_0^i$ is a lattice translate of $\D_0$ for either $i=0$ or $i=1$.  As was previously noted, this must be the case if $X$ is smooth.  Let $\mathfrak{U}=\{\tv(\sigma)\ |\ \sigma\in\Sigma^{(n)}\}$; one can check that this is a subcover of the open affine covering $\mathfrak{U}$ considered in the previous section. We can define then vectors $\lambda_\sigma\in N$ and integers $a_\sigma\in \{-1,1\}$ exactly as in the previous section. Furthermore, set $\lambda_\sigma'=(\lambda_\sigma,\frac{1}{2})$.

\begin{thm}\label{thm:toricks}
	The deformation $\pi$ is locally trivial and its image under the Kodaira-Spencer map is the cocycle defined by
	$$d_{\sigma,\tau}=\sum_{i=1}^n \langle a_\sigma\lambda_\sigma'-a_\tau\lambda_\tau',e_i^*\rangle \chi^{e_i^*-e_n^*} \frac{\partial}{\partial \chi^{e_i^*}}. $$
	In particular, $\pi$ is homogeneous of degree $-e_n^*$.
	\begin{proof}
		The assumption that a polyhedral divisor $\D$ with complete locus has at most non-trivial coefficients for $P=0$ and $\infty$ automatically holds for any toric variety, so we can apply theorem~\ref{thm:tvarks}. The formula follows then directly.
	\end{proof}
\end{thm}

\begin{ex}
We revisit the example from the beginning of this section. Restricting in $\base$ to either $t_{\infty,1}=0$ or $t_{0,1}=0$ gives one-parameter deformations $\pi_0$ or $\pi_\infty$ corresponding, respectively, to the decomposition  $\dfan_0=\dfan_0^0+\dfan_0^1$ or  $\dfan_\infty=\dfan_\infty^0+\dfan_\infty^1$. For the first decomposition, we have $a_\B=-1$, $a_\A=a_\C=a_\D=\ldots=a_\mcH=1$ and $\lambda_\A=-1$, $\lambda_\B=\lambda_\C=\ldots=\lambda_\mcH=0$.
Using the above theorem, we have that the image of $\pi_0$ is described by the cocycle $d$, where
\begin{align*}
d_{\A,\B}=-\chi^{e_1^*-e_2^*}\frac{\partial}{\partial \chi^{e_1^*}}+\frac{\partial}{\partial \chi^{e_2^*}}\qquad d_{\A,\diamond}=-\chi^{e_1^*-e_2^*}\frac{\partial}{\partial \chi^{e_1^*}}\qquad d_{\B,\diamond}=-\frac{\partial}{\partial \chi^{e_2^*}}
\end{align*}
for any $\diamond\neq\A,\B$, and all other terms vanishing.
Likewise, the image of $\pi_\infty$ is described by the cocycle $d$, where
\begin{align*}
d_{\E,\F}=-\chi^{e_1^*+e_2^*}\frac{\partial}{\partial \chi^{e_1^*}}+\frac{\partial}{\partial \chi^{-e_2^*}}\qquad d_{\E,\diamond}=-\chi^{e_1^*+e_2^*}\frac{\partial}{\partial \chi^{e_1^*}}\qquad d_{\F,\diamond}=-\frac{\partial}{\partial \chi^{-e_2^*}}
\end{align*}
for any $\diamond\neq\E,\F$, and all other terms vanishing.
 Thus, the deformation $\pi$ combines deformations of degree $-e_2^*$ and $e_2^*$.

Now, as mentioned in section~\ref{sec:affinetoric}, $X$ is the minimal resolution of a toric Fano surface $X'$ whose fan has rays through $(1,1)$, $(1,-1)$, $(-1,1)$, and $(-1,-1)$. We can then blow down the deformation $\pi$ to a deformation $\pi'$ of $X'$. One can check that this is in fact the same deformation $\pi'$ constructed in the example in section~\ref{sec:affinetoric} by deforming the cone over $X'$ and descending to the quotient.
\end{ex}

Assume now additionally that $X=\tv(\Sigma)$ is smooth and complete. Our goal is to construct $T$-deformations spanning the space of first-order deformations $T_X^1$, which in this case is equal to $H^1(X,\T_X)$. Choose some $R\in M'$  and let  $\rho\in\Sigma^{(1)}$ be some ray with $\langle \rho,R\rangle=1$. Note that by abuse of notation we denote a ray and its primitive generator by the same symbol.
By proper choice of the basis $e_1\ldots,e_n$ above, we can assume that $R=e_n^*$ and that $\rho=e_n$.

Let $\Gamma_{\rho}(-R)$ be the graph embedded in $N_\QQ'$ with vertices consisting of primitive lattice generators of rays $\tau\in\Sigma^{(1)}\setminus \rho$ fulfilling $\langle \tau,R\rangle>0$; two vertices $\tau_1$ and $\tau_2$ are connected by an edge if they generate a cone in $\Sigma$. Note that by rescaling with $\RR_{>0}$ we can consider $\Gamma_{\rho}(-R)$ to be embedded in the slice $\dfan_0$, with vertices of $\Gamma_{\rho}(-R)$ corresponding to non-zero vertices of $\dfan_0$ and with two vertices connected by an edge if they are in fact connected by a line segment in 
$\dfan_0$.

Now for $R\in M'$, define $$\Omega(-R)=\big\{\rho\in\Sigma^{(1)}\ \big|\ \langle \rho,R\rangle=1\ \mathrm{and}\ \Gamma_\rho(-R)\neq \emptyset\big\}.$$
Assume that $\rho\in\Omega(-R)$ and choose now some connected component $C$ of $\Gamma_{\rho}(-R)$. This leads to a decomposition of the slice $\dfan_0$ as follows. 
Consider $\Delta\in \dfan_0$. If $\Delta\cap C=\emptyset$, then set $\Delta^{0}=\Delta$ and $\Delta^{1}=\tail(\Delta)$. If instead the intersection $\Delta\cap C$ is nonempty, set $\Delta^{0}=\tail(\Delta)$ and $\Delta^{1}=\Delta$.

\begin{prop}
	The decompositions $\Delta=\Delta^0+\Delta^1$ form an admissible one-parameter Minkowski decomposition of the slice $\dfan_0$.
\end{prop}
	\begin{proof}
Consider $\Delta\in\dfan$.  Note that the decomposition $\Delta=\Delta^0+\Delta^1$ is always admissible.
We now note that if $\Delta$ is nontrivial, then there is exactly one connected component $C(\Delta)$ of $\Gamma_\rho(-R)$ such that $\Delta\cap C(\Delta)\neq \emptyset$. Indeed, due to nontriviality, $\Delta$ has a vertex corresponding to a vertex $\tau\in\Gamma_\rho(-R)$ with $\tau\neq\rho$. Thus, if $C(\Delta)$ is the connected component of $\tau$, $\Delta\cap C(\Delta)\neq \emptyset$. If for some other connected component $C'$ we have $\Delta\cap C'\neq \emptyset$, then some vertex $\tau'$ must lie in this intersection. Then $\tau$ and $\tau'$ would be rays of some common simplicial cone in the fan $\Sigma$ and thus connected in $\Gamma_\rho(-R)$, so we must have $C'=C(\Delta)$.

For a connected component  $C_i$ of $\Gamma_\rho(-R)$, let $\dfan_0(C_i)=\{\Delta \in \dfan_0\mid\Delta \cap C_i\neq \emptyset\}$. Thus, $\dfan_0$ is the disjoint union
$$\dfan_0= \left (\dot{\bigcup}_i \dfan_0(C_i)\right)\dot{\cup} \{\delta \in \dfan_0\mid\delta=\tail(\delta)\}.$$
Now, $\dfan_0^0$ is obtained from $\dfan_0$ by replacing $\dfan_0(C)$ with $\{\tail(\Delta)\mid\Delta\in \dfan_0(C)\}$ and likewise, we obtain $\dfan_0^1$ by replacing $\dfan_0(C_i)$ with $\{\tail(\Delta)\mid\Delta\in \dfan_0(C_i)\}$ for all components $C_i\neq C$.

We claim that $\dfan_0^0$ is a polyhedral complex; the argument for $\dfan_0^1$ is similar. Indeed, both $\dfan_0\setminus \dfan_0(C)$ and $\{\tail(\Delta)\mid\Delta\in \dfan_0(C)\}$ are polyhedral complexes. Also, the boundary faces of $\{\tail(\Delta)\mid\Delta\in \dfan_0(C)\}$ are exactly the boundary faces of $\dfan_0(C)$. Indeed, the boundary faces of $\dfan_0(C)$ clearly are included among the boundary faces of $\{\tail(\Delta)\mid\Delta\in \dfan_0(C)\}$, but these faces already cover the boundary of $\{\tail(\Delta)\mid\Delta\in \dfan_0(C)\}$. The claim that $\dfan_0^0$ is a polyhedral complex then follows readily. This implies property \ref{item:compface} of definition \ref{def:mdcomp}. Property \ref{item:compatible} of definition \ref{def:mdcomp} follows immediately from the construction.
	\end{proof}

Now let $\pi(C,\rho,R)$ be the one-parameter deformation associated to the above decomposition. We can now formulate one of our main results:

\begin{thm}\label{thm:span}
	Let $X$ be a smooth complete toric variety and $T_X^1(-R)$ the space of first-order deformations in degree $-R$ for some $R\in M'$. Then the one-parameter deformations $\pi(C,\rho,R)$ span $T_X^1(-R)$, where $\rho$ ranges over all rays $\rho\in\Omega(-R)$ and $C$ ranges over all connected components of the graphs $\Gamma_\rho(-R)$.
\begin{proof}

	To prove the theorem, we simply calculate the Kodaira-Spencer map for the above deformations and then use the description of $T_X^1(-R)$ from~\cite{ilten09a}.
	For $\rho\in\Omega(-R)$, let $\partial(R,\rho)$ be the derivation taking $\chi^v\mapsto\langle \rho,v \rangle\chi^{v-R}$. If we choose the basis $e_1,\ldots,e_n$ such that $e_n=\rho$ and $R=e_n^*$, then $\partial{(R,\rho)}=\frac{\partial}{\partial \chi^{e_n^*}}$.	
	Applying theorem~\ref{thm:toricks} we then have that the image of $\pi(C,\rho,R)$ is given by $$d_{\sigma,\tau}=\frac{a'_\sigma-a'_\tau}{2}\,\partial(R,\rho)$$ 	where $a_\sigma'=1$ if $\D_0^\sigma\cap C=\emptyset$ and $a_\sigma'=-1$ otherwise. Indeed,
	 it follows from the above construction that $a_\sigma=a_\sigma'$. Furthermore, $\lambda_\sigma=0$ for all $\sigma\in\Sigma^{(n)}$.

	 Now let $f\in H^0(\Gamma_\rho(-R),\mathbb{C})$, where $H^0(\Gamma_\rho(-R),\mathbb{C})$ is the group of global sections of the locally constant sheaf $\CC$ on the embedded graph $\Gamma_\rho(-R)$. Consider $\sigma\in\Sigma^{(n)}$. If $\Gamma_\rho(-R)\cap \sigma=\emptyset$, set $f_\sigma=1$, otherwise set $f_\sigma=f(v)$ for any $v\in \Gamma_\rho(-R)\cap \sigma$. 	
	From the proof of proposition 2.1 in \cite{ilten09a} together with the (dualized) generalized Euler sequence, see \cite{cox:10a} theorem 8.1.6, we then have the exact sequence
\begin{align*}
\begin{CD}
	0@>>>\displaystyle{\bigoplus_{\rho\in\Omega(-R)}}\mathbb{C}@>>>\displaystyle{\bigoplus_{\rho\in\Omega(-R)}}H^0(\Gamma_\rho(-R),\mathbb{C})@>\Phi>>H^1(X,\T_X)(-R)@>>>0
\end{CD}
\end{align*}
where $\Phi$ maps $f\in H^0(\Gamma_\rho(-R),\mathbb{C})$ to the \v{C}ech cocycle $f_{\sigma,\tau}=\frac{1}{2}(f_\sigma-f_\tau)\partial(R,\rho)$. Now, for $\rho\in\Omega(-R)$ and any connected component $C$ in $\Gamma_\rho(-R)$ let $f(C,\rho,R)\in H^0(\Gamma_\rho(-R))$ be defined by $f(C,\rho,R)_{|C}\equiv -1$ and $f(C,\rho,R)_{|\Gamma_\rho(-R)\setminus C}\equiv 1$. Then we have that $\Phi(f(C,\rho,R))$ is equal to the image of $\pi(C,\rho,R)$ in $H^1(X,T_X^1)$ by the above calculation. Furthermore,  one easily sees that the $f(C,\rho,R)$ form a basis of $H^0(\Gamma_\rho(-R),\mathbb{C})/\mathbb{C}$, where $C$ ranges over all connected components of $\Gamma_\rho(-R)$ except one. Thus, if we allow $\rho$ to vary over the elements of $\Omega(-R)$ as well, the $\pi(C,\rho,R)$ span $T_X^1(-R)$.
\end{proof}
\end{thm}

\remark{Using essentially the same proof, the above result can be extended to the case where $X=\tv(\Sigma)$ for $\Sigma$ a smooth fan with convex support of full dimension, see \cite{ilten:10a} theorem 4.3.2.  In particular, $T$-deformations span $T_X^1$ when $X$ is an equivariant resolution of a toric singularity with no torus factors.} 

\begin{figure}[htbp]
    \centering
    \subfigure[$\dfan_0$]{\threefoldnull}
    \subfigure[$\dfan_\infty$]{\threefoldinfty}
\caption{Slices of $\dfan$ for a toric threefold}\label{fig:3dex}
\end{figure}

\begin{figure}[htbp]
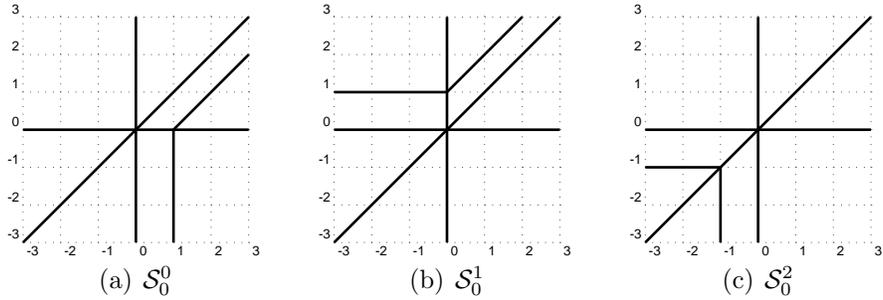

    \centering
    \subfigure[$\dfan_0^0$]{\threefoldnulla}
    \subfigure[$\dfan_0^1$]{\threefoldnullb}
    \subfigure[$\dfan_0^2$]{\threefoldnullc}
\caption{The versal deformation of a toric threefold}\label{fig:3dexvers}
\end{figure}

\begin{ex}
	We shall now construct the versal deformation of a certain toric threefold $X$. Let $X$ be represented as a $T$-variety via the divisorial fan $\dfan$ whose slices are shown in figure~\ref{fig:3dex}, where polyhedra in $\dfan_0$ and $\dfan_\infty$ with equal full-dimensional tailcones belong to a single polyhedral divisor and full-dimensional polyhedra in $\dfan_0$ with one-dimensional tailcone belong to a polyhedral divisor with coefficient $\emptyset$ at $\infty$. If we require that all $\D\in\dfan$ have trivial coefficients everywhere except possibly $0$ and $\infty$, the above uniquely determines $\dfan$ as a divisorial fan.

  We can also describe $X$ in toric terms via a fan $\Sigma$: Let $\rho_1=(1,0,1)$, $\rho_2=(1,1,0)$, $\rho_3=(0,1,1)$, $\rho_4=(-1,0,0)$, $\rho_5(-1,-1,1)$, $\rho_6=\rho_0=(0,-1,0)$, $\rho_7=(0,0,1)$, and $\rho_8=-\rho_7$. The fan $\Sigma$ has top-dimensional cones generated by $\rho_i,\rho_{i+1},\rho_7$ or by $\rho_i,\rho_{i+1},\rho_8$ for $0\leq i <6$. 
	
	In~\cite{ilten09a}, the first author stated that $\dim T_X^1=\dim T_X^1(-R)=2$ for $R=[0,0,1]$. Now, consider the two-parameter deformation $\pi$ corresponding to the decomposition of $\dfan_0=\dfan_0^0+\dfan_0^1+\dfan_0^2$ in figure~\ref{fig:3dexvers}. This is in fact the versal deformation of $X$. Indeed, restricting to $t_{0,i}=0$ for $i=1,2$ gives deformations $\pi(C(i),\rho_7,R)$, where $C(1)=(-1,-1)$ and $C(2)=(0,1)$ are two of the three connected components of $\Gamma_{\rho_7}(-R)$. Thus, the map $T_{\mathbb{A}^2,0}\to T_X^1$ determined by $\pi$ is surjective and so $\pi$ is versal.
\end{ex}

\bibliography{tvardef}
\address{Nathan Ilten\\
Max-Planck-Institut für Mathematik\\
Vivatsgasse 7\\
53111 Bonn, Germany\\}{nilten@cs.uchicago.edu}

\address{Robert Vollmert\\
Mathematisches Institut\\
Freie Universit\"at Berlin\\
Arnimallee 3\\
14195 Berlin, Germany}{vollmert@math.fu-berlin.de}

\end{document}